\documentclass[reqno,10pt,a4paper,final]{amsart}

\usepackage[a4paper,left=35mm,right=35mm,top=30mm,bottom=30mm,marginpar=25mm]{geometry} 
\usepackage{amsmath}
\usepackage{amssymb}
\usepackage{amsthm}
\usepackage{amscd}
\usepackage[ansinew]{inputenc}
\usepackage{bbm}
\usepackage{color}
\usepackage[final]{graphicx}

\renewcommand{\epsilon}{\varepsilon}
\renewcommand{\phi}{\varphi}

\numberwithin{equation}{section}

\newtheoremstyle{thmlemcorr}{10pt}{10pt}{\itshape}{}{\bfseries}{.}{10pt}{{\thmname{#1}\thmnumber{ #2}\thmnote{ (#3)}}}
\newtheoremstyle{thmlemcorr*}{10pt}{10pt}{\itshape}{}{\bfseries}{.}\newline{{\thmname{#1}\thmnumber{ #2}\thmnote{ (#3)}}}
\newtheoremstyle{defi}{10pt}{10pt}{\itshape}{}{\bfseries}{.}{10pt}{{\thmname{#1}\thmnumber{ #2}\thmnote{ (#3)}}}
\newtheoremstyle{remexample}{10pt}{10pt}{}{}{\bfseries}{.}{10pt}{{\thmname{#1}\thmnumber{ #2}\thmnote{ (#3)}}}
\newtheoremstyle{ass}{10pt}{10pt}{}{}{\bfseries}{.}{10pt}{{\thmname{#1}\thmnumber{ A#2}\thmnote{ (#3)}}}

\theoremstyle{thmlemcorr}
\newtheorem{theorem}{Theorem}
\numberwithin{theorem}{section}
\newtheorem{lemma}[theorem]{Lemma}
\newtheorem{corollary}[theorem]{Corollary}
\newtheorem{proposition}[theorem]{Proposition}

\theoremstyle{thmlemcorr*}
\newtheorem{theorem*}{Theorem}
\newtheorem{lemma*}[theorem]{Lemma}
\newtheorem{corollary*}[theorem]{Corollary}
\newtheorem{proposition*}[theorem]{Proposition}
\newtheorem{problem*}[theorem]{Problem}
\newtheorem{conjecture*}[theorem]{Conjecture}

\theoremstyle{defi}
\newtheorem{definition}[theorem]{Definition}

\theoremstyle{remexample}
\newtheorem{remark}[theorem]{Remark}
\newtheorem{example}[theorem]{Example}

\theoremstyle{ass}
\newtheorem{assumption}{Assumption}

\newcommand{\Wrm}{\mathrm{W}}

\newcommand{\Acal}{\mathcal{A}}

\newcommand{\Fcal}{\mathcal{F}}

\newcommand{\Hcal}{\mathcal{H}}

\newcommand{\Lcal}{\mathcal{L}}
\newcommand{\Mcal}{\mathcal{M}}

\newcommand{\Pcal}{\mathcal{P}}
\newcommand{\Qcal}{\mathcal{Q}}

\newcommand{\Scal}{\mathcal{S}}

\newcommand{\Ucal}{\mathcal{U}}

\newcommand{\Cfrak}{\mathfrak{C}}

\newcommand{\Lfrak}{\mathfrak{L}}

\newcommand{\Abb}{\mathbb{A}}

\newcommand{\Mbb}{\mathbb{M}}

\newcommand{\Pbb}{\mathbb{P}}
\newcommand{\Qbb}{\mathbb{Q}}

\newcommand{\Sbb}{\mathbb{S}}

\DeclareMathOperator{\id}{id}

\DeclareMathOperator{\range}{range}

\DeclareMathOperator{\diverg}{div}
\DeclareMathOperator{\curl}{curl}
\DeclareMathOperator{\dist}{dist}
\DeclareMathOperator{\rank}{rank}

\DeclareMathOperator{\supp}{supp}

\newcommand{\ee}{\mathrm{e}}
\newcommand{\ii}{\mathrm{i}}

\newcommand{\set}[2]{\left\{\, #1 \ \textup{\textbf{:}}\ #2 \,\right\}}

\newcommand{\setb}[2]{\bigl\{\, #1 \  \textup{\textbf{:}}\  #2 \,\bigr\}}

\newcommand{\norm}[1]{\|#1\|}

\newcommand{\normb}[1]{\bigl\|#1\bigr\|}

\newcommand{\abs}[1]{|#1|}

\newcommand{\absn}[1]{|#1|}
\newcommand{\absb}[1]{\bigl|#1\bigr|}
\newcommand{\absB}[1]{\Bigl|#1\Bigr|}
\newcommand{\absBB}[1]{\biggl|#1\biggr|}

\newcommand{\dpr}[1]{\langle #1 \rangle}

\newcommand{\dd}{\;\mathrm{d}}

\newcommand{\N}{\mathbb{N}}
\newcommand{\R}{\mathbb{R}}

\newcommand{\Z}{\mathbb{Z}}
\newcommand{\T}{\mathbb{T}}

\newcommand{\rel}{\mathrm{rel}}

\newcommand{\weakly}{\rightharpoonup}
\newcommand{\weaklystar}{\overset{*}\rightharpoondown}

\newcommand{\todown}{\downarrow}
\newcommand{\Lin}{\mathrm{Lin}}
\newcommand{\eps}{\epsilon}

\newcommand{\conv}{\star}

\newcommand{\sbullet}{\begin{picture}(1,1)(-0.5,-2)\circle*{2}\end{picture}}
\newcommand{\frarg}{\,\sbullet\,}
\newcommand{\toYM}{\overset{Y}{\to}}

\newcommand{\proofstep}[1]{\textit{#1}}

\def\Xint#1{\mathchoice 
{\XXint\displaystyle\textstyle{#1}}%
{\XXint\textstyle\scriptstyle{#1}}%
{\XXint\scriptstyle\scriptscriptstyle{#1}}%
{\XXint\scriptscriptstyle\scriptscriptstyle{#1}}%
\!\int} 
\def\XXint#1#2#3{{\setbox0=\hbox{$#1{#2#3}{\int}$} 
\vcenter{\hbox{$#2#3$}}\kern-.5\wd0}} 
\def\dashint{\,\Xint-}

\newcommand{\assref}[1]{A\ref{#1}}

\title{Thin-film limits of functionals on $\mathcal{A}$-free vector fields}

\author{Carolin Kreisbeck}
\address{Carolin Kreisbeck: Fakult\"{a}t f\"{u}r Mathematik, Universit\"{a}t Regensburg, 93040 Regensburg, Germany.}
\email{Carolin.Kreisbeck@mathematik.uni-regensburg.de}

\author{Filip Rindler}
\address{Filip Rindler: Mathematics Institute, University of Warwick, Coventry CV4 7AL, United Kingdom.}
\email{F.Rindler@warwick.ac.uk}

\begin{document}


\begin{abstract}
This paper deals with variational principles on thin films subject to linear PDE constraints represented by a constant-rank operator $\Acal$. We study the effective behavior of integral functionals as the thickness of 
the domain tends to zero, investigating both upper and lower bounds for the $\Gamma$-limit. Under certain conditions 
we show that the limit is an integral functional and give an explicit formula. The limit functional turns out to be 
constrained to $\Acal_0$-free vector fields, where the limit operator $\Acal_0$ is in general not of constant rank.
This result extends work by Bouchitt\'e, Fonseca and Mascarenhas~[\textit{J.\ Convex Anal.} 16 (2009), pp.~351--365] 
to the setting of $\Acal$-free vector fields.
While the lower bound follows from a Young measure approach together with a new
decomposition lemma, the construction of a recovery sequence relies on algebraic considerations
in Fourier space. This part of the argument requires a careful analysis of the limiting behavior of the rescaled 
operators $\Acal_\eps$ by a suitable convergence of their symbols, as well as an explicit construction for plane 
waves inspired by the bending moment formulas in the theory of (linear) elasticity. 
We also give a few applications to common operators $\Acal$.

\vspace{8pt}

\noindent\textsc{MSC (2010):} 49J45 (primary); 35E99, 74K35.
 
\noindent\textsc{Keywords:} dimension reduction, thin films, PDE constraints, $\Acal$-quasiconvexity, $\Gamma$-con\-ver\-gence.

\vspace{8pt}

\noindent\textsc{Date:} \today.
\end{abstract}

\maketitle

\section{Introduction}
For a bounded Lipschitz domain $\omega \subset \R^{d-1}$ and a given (small) thickness $\eps > 0$, define 
$\Omega_\eps := \omega \times (0,\eps)$. The aim of this work is to examine the thin-film limit as $\eps \todown 0$ for 
the variational principles
\begin{equation} \label{eq:var_princ}
  G_\eps[v] \to \min,  \qquad v \colon \Omega_\eps \to \R^m \text{ with $\Acal v=0$ in $\Omega_\eps$,}
\end{equation}
where the functionals $G_\eps$ take the form
\[
  G_\eps[v] = \frac{1}{\eps} \int_{\Omega_\eps} g(y',v(y))\dd{y}.
\]
Here, $y \in \R^d$ is split as $y = (y',y_d)$ and $\Acal$ is the linear first order partial differential operator
\[
  \Acal v := \sum_{k=1}^d A^{(k)} \partial_k v
  \qquad\text{with $A^{(1)}, \ldots, A^{(d)} \in \R^{l\times m}$.}
\]
In applications, $\Omega_\eps$ corresponds to the reference configuration of a thin film with thickness $\eps > 0$, 
the functional $G_\eps$ models the energy stored in a given vector field $v \colon \Omega_\eps \to \R^m$ 
(e.g.\ a deformation or a magnetic field), and the PDE constraint $\Acal v = 0$ encapsulates conditions for admissible 
vector fields $v$. For example, gradients are characterized using $\Acal = \curl$ (as long as $\omega$ is simply connected), 
whereas for solenoidal (incompressible) fields we employ $\Acal = \diverg$. Characterizing the $\Gamma$-limit 
of the problems above then 
corresponds to identifying the effective physical behavior of the system when the thickness $\eps$ goes to zero. Notice 
that the energy density $g$ is assumed not to depend on $y_d$; this reflects the modeling assumption that the film is homogeneous 
with respect to the thickness variable. In fact, it is also possible to treat a dependency on $y_d/\eps$, see below.

The variational treatment of dimension reduction for functionals depending on gradients was initiated by Le~Dret and 
Raoult~\cite{LR95,LR96,LR00}, who rigorously derived the theory of elastic membranes from a three-dimensional elastic model. 
Since then there have been many contributions in the mathematical literature, for 
example~\cite{BFM03, BFM09, BFF00, FF01, Shu00, FJM06}. They include the treatment of different scalings, non-flat limiting 
surfaces, and inhomogeneous materials. In~\cite{BFM09} Bouchitt\'e, Fonseca and Mascarenhas studied a model in elasticity theory that incorporates bending 
by keeping track of deformation away from the mid-plane  in the form of a Cosserat vector, which then appears as an internal 
variable in the energy functional. This work gives a representation formula for the 
corresponding thin-film $\Gamma$-limit. 

Dimension reduction for solenoidal vector fields was recently discussed in~\cite{Kro12}. In this situation, 
the $\Gamma$-limit of the associated energy turns out to be a local functional, characterized entirely by the 
convexification of the energy density. This is in principle due to the fact that the constraint $\diverg v=0$ 
is too weak to prevent the formation of arbitrary oscillations. 
Recent work in the context of Ginzburg--Landau-type functionals can be found in~\cite{ABG10,CS10}. 
In \cite{Kre_micromag}, some of the techniques developed 
below were used to provide an alternative approach to the treatment of thin films in micromagnetics as proposed in \cite{GJ97}.

The idea to work with general PDE constraints can be traced back to the theory of compensated compactness 
introduced by Tartar and Murat~\cite{Tar79,Mur78,Mur81}. The variational theory seems to have started with 
Dacorogna's article~\cite{Dac82} and was further developed by Fonseca and M\"uller~\cite{FM99}, who extensively 
investigated the issue of lower semicontinuity for functionals with $\Acal$-quasiconvex integrands. 
Working in such a general framework allows one to consider a variety of questions, e.g.\ in continuum mechanics 
and electromagnetism (or even both at the same time), in a unified way. Some problems that have already been 
treated within this framework, including relaxation and homogenization, can be found in~\cite{FM99, BFL00, FLM04, FK09}.

\begin{figure}
\label{fig:domain_eps_1}
\def\svgwidth{\textwidth}
\begingroup
  \makeatletter
  \providecommand\color[2][]{%
    \errmessage{(Inkscape) Color is used for the text in Inkscape, but the package 'color.sty' is not loaded}
    \renewcommand\color[2][]{}%
  }
  \providecommand\transparent[1]{%
    \errmessage{(Inkscape) Transparency is used (non-zero) for the text in Inkscape, but the package 'transparent.sty' is not loaded}
    \renewcommand\transparent[1]{}%
  }
  \providecommand\rotatebox[2]{#2}
  \ifx\svgwidth\undefined
    \setlength{\unitlength}{493.76704102pt}
  \else
    \setlength{\unitlength}{\svgwidth}
  \fi
  \global\let\svgwidth\undefined
  \makeatother
  \begin{picture}(1,0.42257289)%
    \put(0,0){\includegraphics[width=\unitlength]{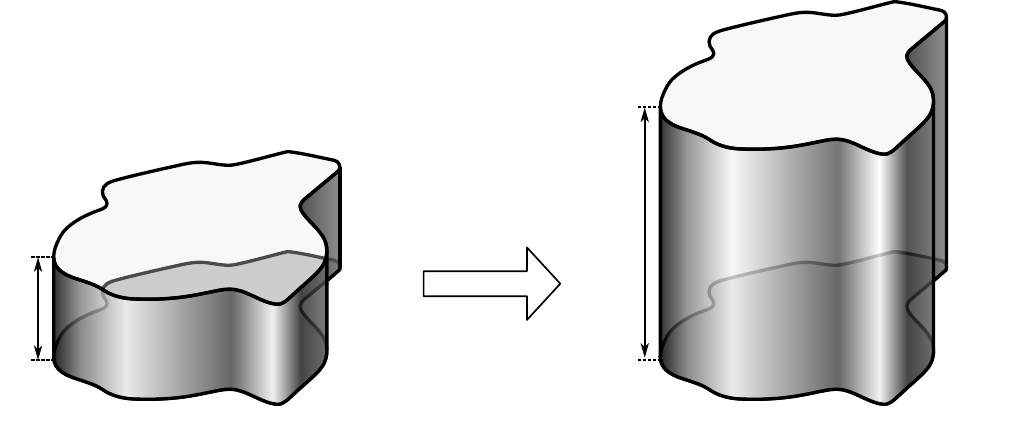}}%
    \put(-0.00134489,0.11686716){\color[rgb]{0,0,0}\makebox(0,0)[lb]{\smash{$\epsilon$}}}%
    \put(0.58840689,0.19787702){\color[rgb]{0,0,0}\makebox(0,0)[lb]{\smash{$1$}}}%
    \put(0.17319515,0.1931659){\color[rgb]{0,0,0}\makebox(0,0)[lb]{\smash{$\omega$}}}%
    \put(0.30834505,0.03307728){\color[rgb]{0,0,0}\makebox(0,0)[lb]{\smash{$\Omega_\epsilon$}}}%
    \put(0.76294695,0.33898365){\color[rgb]{0,0,0}\makebox(0,0)[lb]{\smash{$\omega$}}}%
    \put(0.89809685,0.03307728){\color[rgb]{0,0,0}\makebox(0,0)[lb]{\smash{$\Omega_1$}}}%
    \put(0.13175003,0.00382107){\color[rgb]{0,0,0}\makebox(0,0)[lb]{\smash{$\Acal v = 0$}}}%
    \put(0.72150183,0.00382107){\color[rgb]{0,0,0}\makebox(0,0)[lb]{\smash{$\Acal_\epsilon u = 0$}}}%
  \end{picture}%
\endgroup
\caption{Transformation of domains.}
\end{figure}
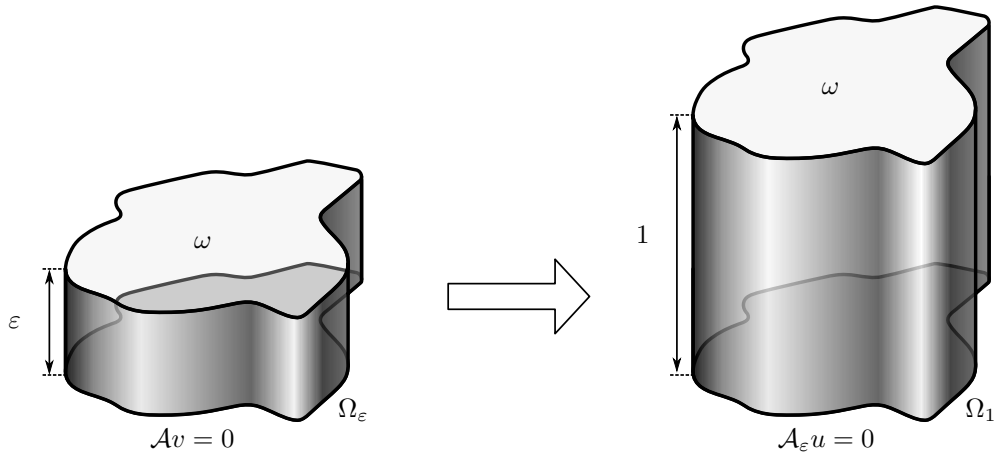

To precisely state our main result, we first transform~\eqref{eq:var_princ} into a problem on the fixed domain 
$\Omega_1=\omega \times (0,1)$ (see Figure~\ref{fig:domain_eps_1}) via the parameter transformation
\begin{equation}\label{rescaling_variables}
  y = (y',y_d)=(x', \eps x_d)  \qquad\text{and}\qquad  u(x) = v(y) = v(x',\eps x_d).
\end{equation}
This re-scaling transforms the PDE constraint $\Acal v = 0$ into
\[
  \Acal_\eps u = 0,  \qquad\text{where}\qquad
  \Acal_\eps u := \sum_{k=1}^{d-1} A^{(k)} \partial_k u + \frac{1}{\eps} A^{(d)}\partial_d u.
\]
Let $f \colon \Omega_1 \times \R^m \to \R$ be a Carath\'{e}odory integrand that satisfies the $p$-growth and $p$-coercivity
conditions
\begin{equation}\label{f_growth/coercivity}
  c \abs{v}^p - C\leq f(x,v) \leq C(1+\abs{v}^p) \qquad \text{for almost every~$x \in \Omega_1$ and all $v \in \R^m$,}
\end{equation}
where $p \in (1,\infty)$ and $c, C > 0$ are constants. We consider the functionals
\begin{align*}
  F_\eps[u] &= \begin{cases}
    \displaystyle\int_{\Omega_1} f(x,u(x))\dd{x}  &\text{if } u \in \Ucal_\eps, \\
    +\infty &\text{otherwise,}
  \end{cases} \\
  \Ucal_\eps &= \setb{ u\in L^p(\Omega_1;\R^m) }{ \text{$\Acal_\eps u=0$ in $\Omega_1$ (distributionally)} }.
\end{align*}
Functions in $\Ucal_\eps$ are called ``$\Acal_\eps$-free'', see Section~\ref{sec:operator_A} for more details. 
Observe that we permit $f$ to depend on $x_d$; this case corresponds to $g$ in~\eqref{eq:var_princ} depending 
on $y_d/\eps$; see~\cite{Shu00} for an application of this. If $f(x,\frarg) = f(x',\frarg)$, however, we have $f = g$.

The main task now is to establish $\Gamma$-upper and lower bounds for $F_\eps$ as $\eps \todown 0$ with respect to the weak 
topology in $L^p(\Omega_1;\R^m)$ and, in some cases, calculate the $\Gamma$-limit. First, we need to define a suitable 
limit operator $\Acal_0$ such that if $u_\eps \weakly u_0$ in $L^p(\Omega_1;\R^m)$ as $\eps \todown 0$ and $\Acal_\eps u_\eps = 0$, 
then $\Acal_0 u = 0$. With the notation $M^i$ (or $[M]^i$) for the $i$th row of the matrix $M$, the operator $\Acal_0$ turns out to be
\begin{align}\label{intro:Acal0}
  \Acal_0 u := \left( \left\{ \begin{aligned}
      &[A^{(d)}]^i \partial_d u              && \text{if $[A^{(d)}]^i \neq 0$,} \\
      &\sum_{k=1}^{d-1} [A^{(k)}]^{i} \partial_k u && \text{if $[A^{(d)}]^i = 0$}
    \end{aligned} \right\} \right)^{i = 1,\ldots,l}.
\end{align}
Of course, this definition depends on the form in which $\Acal$ is written, this point will be clarified in Assumption~\assref{ass:Ad_triangle} below.
In the case $\Acal=\diverg$, this yields $\diverg_0 u=\partial_d u_d$, while for $\Acal=\curl:=\nabla \times$ ($d=3$) one finds
\[
  \curl_0 u =0  \qquad\text{if and only if}\qquad
  \left(-\partial_3 u_2, \partial_3 u_1, \partial_1 u_2-\partial_2 u_1\right)=0.
\]
Hence, we observe that the operator $\Acal$ is in general ``lost'' in the limit. This is closely related to the fact that 
while our functionals are originally defined on thinner and thinner domains, $\Ucal_0$ might contain functions $u$ that 
are not $(d-1)$-dimensional in the sense that $\partial_d u$ does not necessarily vanish. In particular, minimizers
of the limit functional do not have to be $(d-1)$-dimensional. While counter-intuitive, this effect is natural from a physical point of view. 
In elasticity for example it relates to the theory of the Cosserat vector; see Section~\ref{sec:application} for further comments.

Let $\Qcal_{\Acal} f$ be the $\Acal$-quasiconvex envelope of $f$ with respect to the second argument, i.e.\ 
\[
\Qcal_{\Acal} f(x,v) = \inf\set{\dashint_{Q^d} f(x,v+ w(y))\dd{y}}
{ w\in C^{\infty}(\T^d;\R^m)\cap \ker_{\T^d}\Acal,  \textstyle\int_{Q^d} w \dd y = 0}
\]
for $v\in \R^m$ and $x\in \Omega_1$, where $\T^d$ denotes the $d$-torus, which results from $Q^d:=(0,1)^d$
by gluing together opposite sides.
We define $\Qcal_{\Acal_0} f$ analogously. 
Moreover, for $v\in \R^m$ and $x\in \Omega_1$ set the asymptotic $\Acal_0$-quasiconvex envelope $\Qcal^\infty_{\Acal_0}$ of $f$ to be
\begin{align}\label{def:Qcal_infty}
  \Qcal^\infty_{\Acal_0} f(x,v) = \lim_{\eta \to \infty} \Qcal^\eta_{\Acal_0} f(x,v) = \sup_{\eta > 0} \Qcal^\eta_{\Acal_0} f(x,v),
\end{align}
where
\begin{align*}
\Qcal^\eta_{\Acal_0} f(x,v) &:= \inf\biggl\{ \dashint_{Q^d} f(x,v+ w(y))\dd{y}\ \textup{\textbf{:}}\ 
 w \in C^{\infty}(\T^d;\R^m),
\\ &\qquad\qquad\qquad\qquad\qquad \eta \norm{\Acal_0 w}_{W^{-1,1}(\T^d;\R^l)}\leq 1,\;\; \textstyle\int_{Q^d} w \dd y = 0 \biggr\}
\end{align*}
with $W^{-1,1}(\T^d;\R^l)$ denoting the dual of $W_0^{\infty}(\T^d;\R^l)$.

In this work we will prove the following: If $f$ satisfies $\Qcal_\Acal f = \Qcal^\infty_{\Acal_0} f$ (for example if $f$ is convex with respect to the second argument or if only $\Qcal_\Acal f \geq \Qcal^\infty_{\Acal_0} f$, see Remark~\ref{rem:envelopes}),
then, under further assumptions on $\Acal$ stated below,
the $F_\eps$ indeed $\Gamma$-converge to the functional
\begin{align}\label{def:F0}
  F_0[u] &=\begin{cases}
    \displaystyle\int_{\Omega_1} \Qcal_{\Acal} f(x,u(x))\dd{x}  &\text{if  $u\in \Ucal_0$,} \\
    +\infty &\text{otherwise,}
  \end{cases} \\
  \Ucal_0 &= \setb{u\in L^p(\Omega_1;\R^m)}{\Acal_0 u=0\text{ in $\Omega_1$ (distributionally)}}. \notag
\end{align} 
For general $f$, we obtain an upper and lower bound on the $\Gamma$-limit $F_0$ of $F_\eps$ involving the 
$\Acal$-quasiconvex envelopes of $f$ and $\Qcal^\infty_{\Acal_0}f$, respectively, see Theorem~\ref{mainresult_Gamma} below. 
In this general situation the question of whether the limit functional $F_0$ can be represented as an integral functional, i.e.~whether 
$F_0$ is local, remains an open problem. 
Assuming, though, that the $\Gamma$-limit $F_0$ of $F_\eps$ is a-priori known to be local, we show in Section~\ref{sc:locality} the optimality of the upper bound, so that $F_0$ is given 
by the formula~\eqref{def:F0} as well.

Let us now explain the assumptions we impose on $\Acal$:
First, we require the constant-rank property 
(Assumption~\assref{ass:crp}), which was  introduced by Murat in~\cite{Mur81}. This decisive property, 
along with its essential implications, is discussed in detail in Sections~\ref{sec:operator_A} and~\ref{subsubsec:Aeps/A0}. 
Second, Assumption~\assref{ass:Ad_triangle} entails that the shape of $\Acal$ is in a specific sense 
non-degenerate. This assumption can always be achieved and hence is no restriction (compare Example~\ref{ex:curl}). 
%
%
Third, we need to be able to (approximately) extend vector fields that are $\Acal_0$-free in $\Omega_1$ to vector 
fields that are $\Acal_0$-free on the $d$-torus $\T^d$, where we assume without loss of generality $\omega \subset\subset Q^{d-1}$.
The precise requirement is stated in Assumption~\assref{ass:approx_extension}. 
Finally, towards the end of this introduction we comment on the antisymmetry condition we impose through 
Assumption~\assref{ass:antisymmetry}.
In Section~\ref{sec:examples} we show that all these conditions are satisfied for a variety of operators $\Acal$.

With all these preparations, the main result reads as follows:

\begin{theorem}\label{mainresult_Gamma} 
Let $\Omega_1=\omega\times (0,1)\subset\R^d$ be an open, bounded Lipschitz domain and let 
$f\colon \Omega_1 \times \R^m \to \R$ be a Carath{\'e}odory function satisfying $\eqref{f_growth/coercivity}$. 
Further, suppose that Assumptions~\assref{ass:crp}\,--\,\assref{ass:antisymmetry} below
hold for $\Acal$. 
Then:
\begin{itemize}
  \item[(i)] If $u_j \in \Ucal_{\eps_j}$ ($j \in \N$) and $u \in  L^p(\Omega_1;\R^m)$ are such that 
$u_j\weakly u$ in $L^p(\Omega_1;\R^m)$, then $u\in \Ucal_0$ and it holds that
\[
\int_{\Omega_1} \Qcal^\infty_{\Acal_0} f\bigl(x,u(x)\bigr)\dd{x} \leq \liminf_{j\to \infty} F_{\eps_j}[u_j].
\]
  \item[(ii)] For every $u\in \Ucal_0$ and $\eps_j\todown 0$ for $j\to \infty$, there exists a sequence $u_j \in \Ucal_{\eps_j}$ ($j \in \N$) 
such that $u_j\weakly u$ in $L^p(\Omega_1;\R^m)$ and 
\[
\limsup_{j\to\infty} F_{\eps_j}[u_j] \leq \int_{\Omega_1} \Qcal_{\Acal} f\bigl(x,u(x)\bigr) \dd{x}. 
\]
\end{itemize}
Moreover, if $\Qcal_\Acal f=\Qcal^\infty_{\Acal_0} f$, 
then $F_\eps$ converges to the functional $F_0$ in the 
sense of $\Gamma$-convergence with respect to the weak topology in $L^p(\Omega_1;\R^m)$. 
\end{theorem}

One relevant special case, where $\Qcal_\Acal f$ and $\Qcal^\infty_{\Acal_0}f$ coincide, is for $\Acal=\diverg$. Indeed, it holds
that $\Qcal_{\diverg} f=\Qcal_{\diverg_0}f= \Qcal^\infty_{\diverg_0}f=f^{\ast\ast}$, where $f^{\ast\ast}$ denotes the convex envelope of $f$.
In this sense, Theorem~\ref{mainresult_Gamma} provides a generalization of~\cite{Kro12}. 
For $\Acal=\curl$ the function $\Qcal_{\curl} f$, which corresponds to the classical quasiconvex envelope, 
is in general strictly larger than
$\Qcal_{\curl_0} f$, which equals the so-called cross-quasiconvex envelope and is an upper bound on $\Qcal^\infty_{\curl_0}f$. 
We discuss the details of the gradient case in
Section~\ref{sec:application}. 

Notice that we did not impose any boundary conditions on $u$. Indeed, even identifying physically meaningful 
conditions turns out to be non-trivial (for example, it might be necessary to require different conditions on 
the parts $\partial \omega \times (0,1)$ and $\omega \times \{0,1\}$). Moreover, only  ``natural'' boundary 
conditions for $\Acal$-free maps will be preserved under (strong or weak) limits. Therefore, in the current 
work we limit ourselves to the situation without boundary conditions.

One cornerstone of our proof of Theorem~\ref{mainresult_Gamma} is a projection result that is formulated on the torus, 
owing to the fact that Fourier series methods play a decisive role in its proof.
Indeed, in Theorem~\ref{theorem_projection} we obtain projection operators $\Pcal_\eps$ onto $\Acal_\eps$-free 
fields by adapting Theorem~2.14 of~\cite{FM99} to a parameter-dependent setting. 

The proof of the lower bound, which can be found in Section~\ref{sec:proof_lower-bound}, employs Young measures.
The key ingredient is a new decomposition lemma for a sequence of (almost) $\Acal_\eps$-free vector fields 
(see~Theorem~\ref{prop:decomposition2}).
This decomposition is then used in a blow-up argument.

The construction of the upper bound in 
Section~\ref{sec:proof_recovery} requires a new technique and 
hinges on algebraic investigations of the symbols of the operators $\Acal, \Acal_\eps$ and $\Acal_0$. This 
approach to construct recovery sequences allows for a quite intuitive reasoning in Fourier space and 
does not seem to have been employed before. In fact, we can explicitly compute the limit of the symbols 
of $\Acal_\eps$ (see Lemma~\ref{lem:P_eps_conv}). It is important to notice that -- contrary to what one 
might expect initially -- this algebraic limit gives rise to a different Fourier multiplier operator than 
the $\Acal_0$ exhibited in $\eqref{intro:Acal0}$. In general it is not of constant rank, and not even a 
constant-coefficient partial differential operator.
However, the two operators only differ for waves in the 
($\R^{d-1} \times \{0\}$)-plane, and for those we can find a recovery sequence by deforming in the remaining dimension. 
This construction is motivated by the bending moment formula in elasticity theory and requires the antisymmetry relation in
Assumption~\assref{ass:antisymmetry}.

Finally, let us remark that while the Fourier methods developed in this paper provide some interesting insights into 
the structure and geometry of the dimension reduction problem, they are also the precise reason why 
(if we want to work on general domains) we need to require the existence of appropriate extension operators in 
the sense of Assumption~\assref{ass:approx_extension}.


\section{Preliminaries and technical tools}\label{sec:tools}


\subsection{Notation}

Let $\omega \subset \R^{d-1}$ be an open, bounded Lipschitz domain and set $\Omega_\eps := \omega \times (0,\eps)$ 
for $\eps>0$. We will always assume without loss of generality that $\omega\subset\subset Q^{d-1}$, where 
$Q^k := (0,1)^k$ denotes the $k$-dimensional open unit cube. Unless stated otherwise, in the following $1<p<\infty$, 
and $p' = p/(p-1)$ is the dual exponent to $p$.

For a matrix $A \in \R^{l \times m}$, we denote by $\abs{A}$ its Frobenius norm, i.e.\ the vector norm on $\R^{lm}$. 
For $x \in \R^d$, let $x'$ be the vector of the first $d-1$ components, $x' = (x_1,\ldots,x_{d-1})$. We employ 
$\mathbf{e}_k$ for the $k$th unit vector in $\R^d$ and designate the unit sphere in $\R^d$ by $\Sbb^{d-1}$. 
Moreover, $\overline{\R}:=\R\cup\{\infty\}$. For the volume of a measurable set $U \subset \R^d$ we use the notation 
$|U|$, meaning $|U|=\Lcal^d(U)$, where $\Lcal^d$ is the $d$-dimensional Lebesgue measure. 

Further, we use the letter $c$ for constants that can be determined from the known quantities. To stress dependence on a 
specific parameter we employ subscripts, for instance $c_p$ indicates that $c$ depends in particular on $p$. Notice that 
the actual values of constants may differ from line to line.

\subsection{Functions on the torus}
We will often work with the $d$-dimensional torus $\T^d$, which is obtained from $Q^d$ by gluing together opposite sides. 
On the torus, we define the space $C(\T^d)$ to contain all functions from the space $C(Q^d)$ that are also continuous 
over the gluing boundaries of $\T^d$. Next, $C^k(\T^d)$ with $k \in \N$ is the space of all $Q^d$-periodic 
$C^k(Q^d)$-functions whose derivatives (up to $k$th order) 
can be continuously extended to $\overline{Q^d}$ and are again $Q^d$-periodic.
The norm on $C^k(\T^d)$ is given by
\begin{equation*}
  \|\phi\|_{C^k(\T^d)} := \sum_{|\alpha| \leq k}
  \max \setb{|\partial^\alpha \phi(x)|}{x\in \overline{Q^d}},
  \qquad  \phi \in C^k(\T^d),
\end{equation*}
where the summation is over all multi-indices $\alpha \in (\N \cup \{0\})^d$ with 
$|\alpha| := \alpha_1 + \ldots + \alpha_d \leq k$. 
As usual, $C^\infty(\T^d)$ is the intersection of all the spaces $C^k(\T^d)$, $k \in \N$. The space $L^p(\T^d)$ for 
$p \in [1,\infty]$ is simply the space $L^p(Q^d)$.

In the following, for technical reasons, we will also use the ``$E$-torus $\T^d(E)$'', which results from identifying opposite sides 
of the open $d$-dimensional cuboid $E\subset\R^d$. The spaces $C^k(\T^d(E))$ with $k\in \N$, $C^{\infty}(\T^d(E))$ and $L^p(\T^d(E))$ 
are defined analogously to the ones on the torus $\T^d$($=\T^d(Q^d)$).

The discrete Fourier coefficients of a function $f\in L^1(\T^d)\cong L^1(Q^d)$ are defined as
\[
\hat{f}(\xi):= \int_{Q^d} f(x)\ee^{-2\pi \ii x\cdot \xi}\dd{x},\qquad \xi \in \Z^d,
\]
where $\Z^d$ in this context is also called the unit lattice.

\subsection{Sobolev spaces on the torus}\label{sec:Sobolev_torus}

For distributions on the torus $\T^d$ the right kind of test functions are the functions in the 
space $C^\infty(\T^d)$; in particular, no compact support condition is imposed. Then, the Sobolev spaces $W^{k,p}(\T^d)$ can be 
defined in two different ways: Assuming $k \in \N \cup \{0\}$, the first option is to define $W^{k,p}(\T^d)$ in the 
integration-by-parts sense with the aforementioned test functions. The second way is through Fourier Analysis, by saying that a 
distribution $u$ on $\T^d$ (an element of the dual space to $C^\infty(\T^d)$) lies in $W^{k,p}(\T^d)$, where now $k \in \R$ is allowed, 
if and only if
\[
  (I-\Delta)^{k/2}u \in L^p(\T^d).
\]
Here $(I-\Delta)^{k/2}$ is the Fourier multiplier operator with symbol $\xi \mapsto (1+4\pi^2\abs{\xi}^2)^{k/2}$. 
Accordingly, the norm on $W^{k,p}(\T^d)$ can be defined either as
\[
  \Biggl( \sum_{\abs{\alpha} \leq k} \norm{\partial^\alpha u}_{L^p(\T^d)}^p \Biggr)^{1/p} \text{ for $k\in \N\cup\{0\}$}
  \qquad\text{or as}\qquad
  \norm{(I-\Delta)^{k/2} u}_{L^p(\T^d)} \quad\text{ for $k\in \R$}.
\]
By virtue of the Mihlin Multiplier Theorem, for $k \in \N \cup \{0\}$, both definitions turn out to be 
equivalent and the corresponding norms are comparable. It can be further shown that for all $k \in \R$ the 
dual space to $W^{k,p}(\T^d)$ is $W^{-k,p'}(\T^d)$ with $1/p + 1/p' = 1$. Again, we have that the two possible 
definitions of the dual norm in $W^{-k,p'}(\T^d)$ for $k \in \N \cup \{0\}$, namely
\[
  \sup_{\phi \in W^{k,p}(\T^d),\phi \neq 0} \frac{\abs{\dpr{u,\phi}}}{\norm{\phi}_{W^{k,p}(\T^d)}}
  \qquad\text{and}\qquad
  \normb{(I-\Delta)^{-k/2} u}_{L^{p'}(\T^d)},
\]
are equivalent. Of course, analogous statements also hold for the vector-valued spaces $W^{k,p}(\T^d;\R^m)$.

In this work, we will use either definition and norm according to the situation at hand. The above assertions are 
standard and proofs can for example proceed along the lines of Chapter~6 in~\cite{Gra09} or Chapter~VI of~\cite{Ste93} 
(which, however, consider Fourier transforms instead of Fourier series).


\subsection{The differential operator $\Acal$ and its symbol}\label{sec:operator_A}

Given matrices $A^{(1)}, \ldots, A^{(d)}\in\R^{l \times m}$, we define a linear partial differential operator of first order
\begin{equation}\label{operator_A}
  \Acal u := \sum_{k=1}^d A^{(k)} \partial_k u,  \qquad u \in C^1(\R^d;\R^m).
\end{equation}
Its symbol is
\[
\Abb(\xi) := \sum_{k=1}^d A^{(k)} \xi_k,  \qquad \xi \in \R^d.
\]

The partial differential operator $\Acal$ can be viewed as a bounded, linear operator 
$\Acal: L^p(\Omega;\R^m)\to  W^{-1,p}(\Omega;\R^l)$ with $\Omega\subset \R^d$ open, if interpreted as
\[
  (\Acal u)[v]:=-\int_{\Omega} u\cdot \Acal^T v\dd{x} 
\]
for all $u\in L^p(\Omega;\R^m)$ and $v\in W_0^{1,p'}(\Omega;\R^l)$, where $\Acal^T:=\sum_{k=1}^d (A^{(k)})^T \partial_k$. 

By the expression ``$\Acal u=0$ in $\Omega$'' we mean
\[
  - \int_{\Omega} u\cdot \Acal^T \phi\dd{x}=0
  \qquad\text{for all $\phi\in C_c^{\infty}(\Omega;\R^l)$.}
\]

Similarly, for functions $u$ on the $d$-torus $\T^d$ the statement ``$\Acal u = 0$ in $\T^d$'' for some $u \in L^p(\T^d;\R^m)$ is 
understood in the sense of distributions on $\T^d$, i.e.\ 
\[
  - \int_{Q^d} u\cdot \Acal^T \phi\dd{x}=0
  \qquad\text{for all $\phi\in C^{\infty}(\T^d;\R^l)$.}
\]
This can be expressed equivalently as the algebraic equations 
\begin{equation*}
  \Abb(\xi) \hat{u}(\xi) = 0  \quad \in \R^l  \qquad\text{for all $\xi \in \Z^d$}.
\end{equation*}

Notice that the condition ``$\Acal u = 0$ in $\T^d$'' 
also includes the requirement that this equation holds over the gluing 
boundaries of $\T^d$ when thinking of $\T^d$ as originating from $Q^d$ by gluing. So the condition ``$\Acal u = 0$ in $Q^d$'' 
(which only uses test functions of the class $C_c^{\infty}(Q^d;\R^l)$) is in general strictly weaker, even if $u$ is 
differentiable in the classical sense.
Analogously to the $\Acal$-freeness in $\T^d$ we define ``$\Acal u = 0$ in $\T^d(E)$'' by duality with test functions in $C^{\infty}(\T^d(E);\R^l)$. 
In all of the following, $\ker_{\Omega} \Acal$ is the set of all $u\in L^1(\Omega;\R^m)$ such that $\Acal u = 0$ in $\Omega$ 
(in the above sense). Analogously, we define $\ker_{\T^d} \Acal$ and $\ker_{\T^d(E)}\Acal$.

As usual in the theory of $\Acal$-free vector fields, we assume the following fundamental condition:

\begin{assumption}[Constant-rank property] \label{ass:crp}
The rank of the matrix $\Abb(\xi)$ is constant for all $\xi \in \R^{d}\setminus\{0\}$, i.e.\ there is $r\in \N$ such that
\[
\rank \Abb(\xi) = r \qquad\text{ for all }\xi \in \R^{d}\setminus\{0\}.
\]
\end{assumption}

We refer to~\cite{Mur81,Tar79,Tar83,FM99} for more information about the constant-rank property and $\Acal$-free vector fields.

\subsection{The operators $\Acal_\eps$ and $\Acal_0$} \label{subsubsec:Aeps/A0}

Before we come to the definition of $\Acal_\eps$ and $\Acal_0$, we state the following technical assumption on $\Acal$, 
which will turn out to be important:

\begin{assumption}[No linearly dependent rows in $A^{(d)}$] \label{ass:Ad_triangle} 
The number of non-zero rows of the matrix $A^{(d)}$ is equal to the rank of $A^{(d)}$.
\end{assumption}

Assuming the constant-rank property Assumption~\assref{ass:crp}, this assumption entails that the number of non-zero rows of 
the matrix $A^{(d)}$ is equal to the rank of $\Abb(\xi)$ for all $\xi \in \R^{d}\setminus\{0\}$ (see Lemma \ref{lemma_rank_Aeps} below).
 
Notice that Assumption~\assref{ass:Ad_triangle} imposes no restriction, because we can always achieve it by Gaussian elimination, 
i.e.\ adding multiples of rows to other rows until Assumption~\assref{ass:Ad_triangle} is satisfied. Depending on the specific 
operator at hand there might also be a more intuitive way of adapting $\Acal$ to Assumption~\assref{ass:Ad_triangle} 
(see Example~\ref{ex:curl} for an explicit discussion of this matter in the case $\Acal=\curl$).

Using the convention that $M^i$ (or $[M]^i$) denotes the $i$th row of the matrix $M$, we define for $\eps > 0$ and $u \in C^1(\R^d;\R^m)$,
\begin{align}\label{Aeps_A0}
  \Acal_\eps u &:= \sum_{k=1}^{d-1} A^{(k)} \partial_k u
    + \frac{1}{\eps} A^{(d)} \partial_d u, \\
  \Acal_0 u &:= \left( \left\{ \begin{aligned}
      &[A^{(d)}]^i \partial_d u              && \text{if $[A^{(d)}]^i \neq 0$,} \\
      &\sum_{k=1}^{d-1} [A^{(k)}]^{i} \partial_k u && \text{if $[A^{(d)}]^i = 0$,}
    \end{aligned} \right\} \right)^{i = 1,\ldots,l}.\nonumber
\end{align}
Notice also that for the line-by-line definition of $\Acal_0$ it is essential that we require the technical assumption above, 
for otherwise we might get a substantially different $\Acal_0$. This can be seen through the following example:
\begin{example}
In two dimensions, consider the constant-rank operator $\Acal$ defined as
\[
  \Acal u := \begin{pmatrix} 0 & 1 \\ 0 & 0 \end{pmatrix} \partial_1 u
  + \begin{pmatrix} 0 & 1 \\ 0 & 1 \end{pmatrix} \partial_2 u.
\]
Clearly, this $\Acal$ does not satisfy Assumption~\assref{ass:Ad_triangle}. If we still apply the previous definition for $\Acal_0$, we get
\[
  \Acal_0 u = \begin{pmatrix} 0 & 1 \\ 0 & 1 \end{pmatrix} \partial_2 u.
\]
This, however, is not the right operator for our purposes: Every smooth function $u = (u_1,u_2)^T$ with $\Acal u = 0$ satisfies not only
\begin{align*}
  \partial_1 u_2 + &\partial_2 u_2 = 0,  \\
  &\partial_2 u_2 = 0,
\end{align*}
but also, by subtracting the second condition from the first,
\[
  \partial_1 u_2 = 0.
\]
For an $\Acal_\eps$-free sequence, this property would clearly also hold for the limit, but in the above definition of 
$\Acal_0$ it does not appear. 
This shows the need to eliminate linearly dependent rows from $A^{(d)}$.
\end{example}

For the symbols of $\Acal_\eps$ and $\Acal_0$, we have
\begin{align*}
  \Abb_\eps(\xi) &:=  \sum_{k=1}^{d-1} A^{(k)} \xi_k + \frac{1}{\eps} A^{(d)} \xi_d   \qquad\text{for $\eps > 0$,}\\
  \Abb_0(\xi) &:= \left( \left\{ \begin{aligned}
      &[A^{(d)}]^i \xi_d              && \text{if $[A^{(d)}]^i \neq 0$,} \\
      &\sum_{k=1}^{d-1} [A^{(k)}]^i \xi_k && \text{if $[A^{(d)}]^i = 0$,}
    \end{aligned} \right\} \right)^{i = 1,\ldots,l}, \qquad \xi \in \R^d.
\end{align*}

Notice that even under Assumption~\assref{ass:crp}, the constant-rank property cannot be guaranteed for $\Acal_0$ 
(see Section~\ref{sec:examples}). Regarding $\Acal_\eps$, however, it is satisfied for all $\eps>0$:

\begin{lemma}\label{lemma_rank_Aeps}
If $\Acal$ is a constant-rank operator in the sense of Assumption~\assref{ass:crp}, then
\[
 \rank \Abb_\eps(\xi)=\rank A^{(d)}
\]
for all $\xi\in \R^d\setminus\{0\}$ and $\eps>0$.
\end{lemma}

\begin{proof}
We define $\xi_\eps:=(\xi',\eps^{-1}\xi_d)$ for $\xi\in \R^d \setminus \{0\}$,
and observe that
\[
 \Abb(\xi_\eps) =  \sum_{k=1}^{d-1} A^{(k)} \xi_k
 + \frac{1}{\eps} A^{(d)} \xi_d  = \Abb_\eps(\xi).
\]
In view of the fact that $\Abb(\mathbf{e}_d) = A^{(d)}$, the constant-rank property of Assumption~\assref{ass:crp} implies
\[
   \rank \Abb_\eps(\xi) = \rank \Abb(\xi_\eps) 
   = \rank \Abb(\mathbf{e}_d)
   = \rank A^{(d)},
\]
from which the assertion follows.
\end{proof}

 
Next let us introduce some notation and provide an auxiliary result in Lemma~\ref{lem:lindep}. Unless stated otherwise we assume in the following that 
$\Acal$ satisfies Assumptions~\assref{ass:crp} and ~\assref{ass:Ad_triangle}. We write
\begin{align}\label{notationA+}
A^{(d)}= \left[\begin{array}{c}A^{(d)}_+ \\ \hline A^{(d)}_-\end{array}\right]
= \left[\begin{array}{c}A^{(d)}_+\\ \hline 0\end{array}\right]
 \qquad\text{and}\qquad
\Abb(\xi)=\left[\begin{array}{c} \Abb(\xi)_+\\ \hline \Abb(\xi)_-\end{array}\right]
= \left[\begin{array}{c}\Abb(\xi)_+\\ \hline \Abb'(\xi')_-\end{array}\right],
\end{align}
where $\xi=(\xi', \xi_d)\in \R^d$, $A^{(d)}_+ \in \R^{r \times m}$, and $\Abb'(\eta) =\sum_{k=1}^{d-1}A^{(k)}\eta_k$ for $\eta\in \R^{d-1}$. 
Notice that $A^{(d)}_+$ has $r$ linearly independent rows,
which are exactly the non-zero rows of $A^{(d)}$ according to Assumption~\assref{ass:Ad_triangle}, and $A^{(d)}_-=0\in \R^{(l-r)\times m}$. 
Then, 
\begin{align*}
\Abb_\eps(\xi)=\left[\begin{array}{c} \Abb_\eps(\xi)_+\\ \hline \Abb_\eps(\xi)_-\end{array}\right]
=\left[\begin{array}{c} \Abb_\eps(\xi)_+ \\ \hline \Abb'(\xi')_-\end{array}\right]
 \qquad\text{and}\qquad
\Abb_0(\xi)=\left[\begin{array}{c} \Abb_0(\xi)_+\\ \hline \Abb_0(\xi)_-\end{array}\right]
=\left[\begin{array}{c} A^{(d)}_+ \xi_d\\ \hline \Abb'(\xi')_-\end{array}\right].
\end{align*}
In the following we use the notation $\Acal_{+}$ and $\Acal_{-}$ to refer to the differential 
operator with symbol $\Abb(\xi)_+$ and $\Abb(\xi)_-$, respectively.

\begin{lemma}\label{lem:lindep}
Let $\Acal$ be a constant-rank operator of the form \eqref{operator_A}.
If $\xi \in\R^d$ with $\xi_d\neq 0$, then $\Abb'(\xi')_-$ is linearly dependent on $A^{(d)}_+$, i.e.\ each row of $\Abb'(\xi')_-$ 
can be written as a linear combination of rows of $A^{(d)}_+$. 
\end{lemma}

\begin{proof}
Since the cases $l=r$ and $m=r$ are trivial, we may assume in the following that $l\geq r+1$ and $m\geq r+1$.
First let us show that, provided $\eps>0$ is sufficiently small, 
\begin{align}\label{lem:step1}
 \rank \Abb_\eps(\xi)_+=r \quad\text{ for $\xi\in \R^d$ such that $\xi_d\neq 0$}. 
\end{align}
We observe that $\eps \Abb_\eps(\xi)_+\to A^{(d)}_+\xi_d$ as $\eps$ tends to $0$. Recall that
$\rank A^{(d)}_+ =r$, select an $(r\times r)$-submatrix $M(A^{(d)}_+)$ with $\det M(A^{(d)}_+) \neq 0$
and let $M(\Abb_\eps(\xi)_+)$ be the corresponding submatrix of $\Abb_\eps(\xi)_+$. Then, 
\begin{align*}
 \eps^r \det M(\Abb_\eps(\xi)_+)=  \det M(\eps \Abb_\eps(\xi)_+) \to \det M(A^{(d)}_+\xi_d)=\xi_d^r \det M(A^{(d)}_+)\neq 0
 \qquad\text{as $\eps\to 0$.}
\end{align*}
Hence, $\rank \Abb_\eps(\xi)_+\geq r$ for $\eps>0$ small enough.
On the other hand, we infer from Lemma~\ref{lemma_rank_Aeps} that
\begin{align*}
 r=\rank \Abb_\eps(\xi)=\rank \left[\begin{array}{c}\Abb_\eps(\xi)_+\\ \hline \Abb'(\xi')_-\end{array}\right],
\end{align*}
which implies $\rank\Abb_\eps(\xi)_+\leq r$ and finishes the proof of \eqref{lem:step1}.

Next we define 
\begin{align*}
  \Mbb_\eps(\xi)=\left[\begin{array}{c} \eps\Abb_\eps(\xi)_+\\ \hline \Abb'(\xi')_-\end{array}\right].
\end{align*}
Then, $\rank \Mbb_\eps(\xi)=
\rank \Abb_\eps(\xi)_+ =r$,
because all rows of $\Abb'(\xi')_-$ can be written as linear combinations of rows of 
$\Abb_\eps(\xi)_+$, and 
\begin{align*}
 \Mbb_\eps(\xi)\to \Abb_0(\xi)\qquad \text{as $\eps\to 0$}.
\end{align*}

Choose any $[(r+1)\times(r+1)]$-submatrix $M(\Abb_0(\xi))$ and $M(\Mbb_\eps(\xi))$ involving the first $r$ rows. 
By the continuity of the determinant, $\det M(\Abb_0(\xi))=\lim_{\eps\to 0}\det M(\Mbb_\eps(\xi))=0$.
The assertion follows, since the
first $r$ rows of $A^{(d)}$ are linearly independent according to Assumption~\assref{ass:Ad_triangle}.
\end{proof}

\begin{remark}
In particular, Lemma~\ref{lem:lindep} implies that $\rank \Abb_0(\xi)=\rank \Abb(\xi) = r$ for all $\xi\in \R^d$ with
$\xi_d\neq 0$.
\end{remark}

\subsection{Approximate extensions and antisymmetry conditions.}

The final two assumptions are necessary for the proof of the upper bound in Section~\ref{sec:proof_recovery}.
As alluded to in the introduction, we have to require that $\Acal_0$-free fields in $\Omega_1$ can be 
extended approximately to be $\Acal_0$-free on the torus. This is necessary to provide the correct context for the application of our
Fourier methods. 

\begin{assumption}[Approximate extension]\label{ass:approx_extension}
Suppose $u\in L^p(\Omega_1;\R^m)$ with $\Acal_0 u = 0$ in $\Omega_1$. Then there exists a sequence 
$(\bar{u}_j)_{j}\subset L^p(Q^d;\R^m)$ such that $\bar{u}_j$ is $\Acal_0$-free in $\T^d$ for each 
$j\in \N$ and $\,\bar{u}_j \to u$ in $L^p(\Omega_1;\R^m)$ as $j\to \infty$.
\end{assumption}
\begin{remark}\label{rem:extension}
For the moment let us allow for more general domains $\Omega \subset \T^d$ (as opposed to our 
standard requirement $\Omega_1 = \omega\times(0,1)$ with $\omega\subset\subset Q^{d-1}$).
Assuming $\Omega = \T^d$, we observe that the extension property is trivially fulfilled for any $\Acal_0$. 
In comparison to employing $\Omega=Q^d$, working on the $d$-torus implicitly imposes additional boundary 
conditions on $\Acal$-free functions. Indeed, for $u\in C^1(\overline{Q^d};\R^m)$ we can integrate by parts to get
\begin{align*}
\int_{Q^d} \Acal u\cdot v\dd{x}= - \int_{Q^d} u\cdot \Acal^T v\dd{x} + \int_{\partial Q^d} \Abb(n)u \cdot v\dd{\Hcal^{d-1}}, 
\qquad v\in C^1(\overline{Q^d};\R^l),
\end{align*}
where $n:\partial Q^d\to \Sbb^{d-1}$ is the exterior unit normal of $Q^d$. Hence, $\Acal u=0$ in $\T^d$ if and only if $\Acal u =0$ 
in $Q^d$ and 
\begin{align*}\Abb(\mathbf{e}_j)u|_{\partial Q^d\cap \{x_j=0\}}
=\Abb(\mathbf{e}_j)u|_{\partial Q^d \cap \{x_j=1\}} \qquad\text{ for all $j=1,\ldots, d$.}
\end{align*}
\end{remark}

Moreover, the construction of a recovery sequence for plane waves requires the following:
\begin{assumption}[Antisymmetry relation]\label{ass:antisymmetry}
The matrices $A^{(k)}\in \R^{l\times m}$ in~\eqref{operator_A}
are such that 
\begin{align}\label{equation:antisymmetry}
 A^{(k)}(A^{(d)})^\dagger A^{(j)} = - A^{(j)} (A^{(d)})^\dagger A^{(k)} \qquad \text{for $k,j=1, \ldots, d-1$,}
\end{align}
where $(A^{(d)})^\dagger\in \R^{m\times l}$ is the Moore-Penrose pseudoinverse of $A^{(d)}$.
\end{assumption}

This assumption in particular entails that $A^{(k)}(A^{(d)})^\dagger A^{(k)} = 0$ for $k=1, \ldots, d-1$.

\subsection{Examples and applications}\label{sec:examples}
In this section we investigate in detail the two most prominent constant-rank operators, which are 
$\diverg$ and $\curl$, and list some more applications in which $\Acal$-free vector fields play a decisive role 
(see also Remark~3.3 of \cite{FM99}, \cite{BFL00}, \cite{Kro12}, \cite{Tar79} and \cite{Tar83}). 
In particular, we discuss Assumptions~\assref{ass:crp}--\assref{ass:antisymmetry}.

\subsubsection{$\Acal=\diverg$}\label{ex:div}
This case corresponds to working on solenoidal vector fields.
Regarding the notation of $\eqref{operator_A}$, here $m=d > 1$ and $l=1$. 
For $u:\R^d\to\R^d$ we have
\begin{align*}
 \diverg u =\nabla\cdot u=\sum_{k=1}^d \partial_k u_k=\sum_{k=1}^d A^{(k)}_{\diverg} \partial_k u
\end{align*}
with $A^{(k)}_{\diverg}:=\mathbf{e}_k^T\in \R^{1\times d}$. Clearly, Assumption~\assref{ass:Ad_triangle} is true. Moreover, $\diverg$ meets Assumption~\assref{ass:antisymmetry}, 
since $A^{(k)}_{\diverg}(A^{(d)}_{\diverg})^\dagger=\mathbf{e}_k^T\mathbf{e}_d=0$ for 
$k=1, \ldots, d-1$. 
The symbol of $\diverg$ reads $\Abb_{\diverg}(\xi)=\xi^T$ with $\xi\in\R^d$, so that
$\ker \Abb_{\diverg}(\xi)=\{v\in \R^d: \xi\cdot v=0\}$. Hence, $\dim \ker \Abb_{\diverg}(\xi)= d-1$ 
for all $\xi\in \R^d\setminus \{0\}$ and $\diverg$ fulfills the constant-rank condition of Assumption~\assref{ass:crp} with $r=1$. 
Now consider for $\eps>0$,
\begin{align*}
 \diverg_\eps u= \sum_{k=1}^{d-1} \partial_k u_k + \frac{1}{\eps} \partial_d u_d
\end{align*}
or in short notation $\diverg_\eps u=\nabla_\eps\cdot u$, 
where $\nabla_\eps=(\partial_1,\ldots, \partial_{d-1}, \eps^{-1}\partial_d)^T$. In view of $\eqref{Aeps_A0}$, 
it holds that $\diverg_0 u=\partial_d u_d$. Then, for $\xi\in \R^d$ one has $\ker \Abb_{\diverg_0} (\xi)=\set{v\in \R^d}{\xi_d v_d=0}$, 
so that
\begin{align*}
\dim \ker \Abb_{\diverg_0}(\xi)=\begin{cases}
                   d-1 & \text{if $\xi_d\neq 0$,}\\
 		   d   & \text{if $\xi_d=0$.}
               \end{cases}
\end{align*}
This entails that $\diverg_0$ is not of constant rank.

Regarding Assumption~\assref{ass:approx_extension} we observe: If $u\in L^p(\Omega_1;\R^m)$ with $\diverg_0 u=\partial_d u_d=0$ in $\Omega_1$, 
then $u_d$ has to be constant with respect to the $x_d$-variable, and extending $u$ by zero to $Q^d$ 
(we call this extension $\bar{u}$) preserves this property.
Thus, $\diverg_0 \bar{u}=0$ in $\T^d$, which means that there is even an exact extension of $u$.

\subsubsection{$\Acal=\curl$.}\label{ex:curl}
The $\curl$ of a matrix-valued function $F:\R^d\to \R^{n\times d}$ is defined row by row as
\begin{align}\label{def_curl}
 (\curl F)_{ijk} =\partial_j F_{k}^i -\partial_k F_{j}^i,
 \qquad  1\leq j,k\leq d \text{ and } 1\leq i\leq n.
\end{align}
Assume that $d>1$.
In terms of $\eqref{operator_A}$ we have $m=n d$, $l=d^2 n$ and 
\begin{align}\label{matrices_Acurl}
\bigl(A_{\curl}^{(r)}\bigr)_{ijk,qp}= \delta_{rj}\delta_{qi}\delta_{pk} - \delta_{rk}\delta_{qi}\delta_{pj}, 
\qquad  1\leq j,k,p,r\leq d \text{ and } 1\leq i,q\leq n.
\end{align}
Moreover,
\begin{align*}
\ker \Abb_{\curl} (\xi)&= \set{G\in \R^{n\times d}}{\xi_j G_{k}^i-\xi_k G_{j}^i = 0, \ 1\leq j,k\leq d \text{ and } 1\leq i\leq n}\\
&=\set{G\in \R^{n\times d}}{G= a\otimes \xi,\ a\in \R^n}, \qquad \xi \in \R^d,
\end{align*}
which shows that $\curl$ meets Assumption~\assref{ass:crp} with $r=n(d-1)$, since $\dim \ker \Abb_{\curl}(\xi)=n$ for all $\xi \in \R^d\setminus \{0\}$. 

As pointed out in Section~\ref{subsubsec:Aeps/A0}, before we can state the correct limit operator $\curl_0$ of $\curl_\eps$, 
we have to check Assumption~\assref{ass:Ad_triangle} first. 
A close look at $\eqref{matrices_Acurl}$ reveals that the number of non-zero rows of $A_{\curl}^{(d)}$ is $2n(d-1)$, 
while the rank of $A_{\curl}^{(d)}$ is only $n(d-1)$. The reason for this is a redundancy in the definition of 
$\curl$ in $\eqref{def_curl}$. To obtain an operator that is equivalent to $\curl$ (in the sense that the kernels coincide) 
and fits into the framework of this work, we need to get rid of this symmetry. Here we simply choose the additional requirement $k<j$ in $\eqref{def_curl}$. Then $l=(d/2)(d-1)n$. In the case $d=3$ we 
can also equivalently use the more natural definition $\curl F^i:= \nabla \times F^i$ (row-wise). Notice that we refer to $\curl$ 
in this new form from now on without change of notation. Then, 
\begin{align*}
\bigl(\curl_\eps F\bigr)_{ijk}=\begin{cases}
\frac{1}{\eps}\partial_d F_{k}^i- \partial_k F_{d}^i & \text{ if $1 \leq k\leq d-1$ and $1\leq i\leq n$}, \\
\partial_j F_{k}^i-\partial_k F_{j}^i & \text{ if $1\leq j,k \leq d-1$ and $k<j$ and $1\leq i\leq n$,}\end{cases}
\end{align*}
and the corresponding limit operator reads
\begin{align}\label{curl0}
\bigl(\curl_0 F\bigr)_{ijk}=\begin{cases}
\partial_d F_{k}^i & \text{ if $1 \leq k \leq d-1$ and $1\leq i\leq n$}, \\
\partial_j F_{k}^i-\partial_k F_{j}^i & \text{ if $1\leq j,k \leq d-1$ and $k<j$ and $1\leq i\leq n$.}\end{cases}
\end{align}

Let us remark in passing that $\curl_0$ is another example of an operator failing to have constant rank. Indeed, for $\xi\in \R^d\setminus\{0\}$, 
\begin{align*}
 \ker \Abb_{\curl_0} (\xi)=\begin{cases}
                            \set{G\in \R^{n\times d}}{G=(0| G_d)} & \text{ if $\xi_d\neq 0$,}\\
                            \set{G\in \R^{n\times d}}{G=(a\otimes \xi'| G_d),\ a\in \R^n} & \text{ if $\xi_d= 0$}.
                           \end{cases}
\end{align*}
Thus, $\dim \ker\Abb_{\curl_0}(\xi)=n$ if $\xi_d\neq 0$ and $\dim \ker \Abb_{\curl_0}(\xi)=2n$ if $\xi_d= 0$.

For $\Acal_0=\curl_0$ one can prove an exact extension result, which implies Assumption~\assref{ass:approx_extension}:

\begin{lemma}\label{lem:extension_curl_free}
Let $\omega \subset\subset Q^{d-1}$ be open, bounded and simply connected, and assume $F\in L^p(\Omega_1;\R^{n\times d})$ is $\curl_0$-free in $\Omega_1$.
Then there exists $\bar{F}\in L^p(\T^d;\R^{n\times d})$ such that $\bar{F}|_{\Omega_1}=F$ and $\curl_0 \bar{F}=0$ in $\T^d$.
\end{lemma}
\begin{proof}
In view of the representation of $\curl_0$ in $\eqref{curl0}$ we find that $F':=(F_1|\cdots|F_{d-1})$ depends only on $x'$ 
and is $\curl$-free in $\omega\subset\R^{d-1}$ (here $\curl$ stands for the $(d-1)$-dimensional $\curl$-operator). Since $\omega$ is simply connected, $F'$ possesses a potential. 
Precisely, there exists $v\in W^{1,p}(\omega;\R^n)$ such that $F'=\nabla v$.
Let $\bar{v}\in W^{1,p}(Q^{d-1};\R^n)$ be a standard Sobolev extension of $v$ with $\supp \bar{v}\subset\subset Q^{d-1}$. 
Then $\bar{v}\in W^{1,p}(\T^{d-1};\R^n)$.
By setting $\bar{F}'=\nabla \bar{v}$ and extending $F_d$ by zero to $Q^d$, we end up with $\bar{F}$,
which is $\curl_0$-free in $\T^d$ by construction.
\end{proof}

Finally, $\curl$ is antisymmetric in the sense of Assumption~\assref{ass:antisymmetry}.
For $d=3$ and $n=1$ with $\curl$ in the form $\curl=\nabla \times$ this is easy to see.
Indeed, in this case one has 
\begin{align*}
A^{(1)}_{\rm curl} = \left(\begin{array}{ccc}0 & 0 & 0\\ 0 & 0 & -1 \\ 0& 1 & 0\end{array}\right), \quad
A^{(2)}_{\rm curl} = \left(\begin{array}{ccc}0 & 0 & 1\\ 0 & 0 & 0 \\ -1& 0 & 0\end{array}\right), \quad
A^{(3)}_{\rm curl} = \left(\begin{array}{ccc}0 & -1 & 0 \\ 1 & 0 & 0 \\ 0 & 0 &0\end{array}\right),
\end{align*}
and $(A^{(3)}_{\rm curl})^\dagger=(A^{(3)}_{\rm curl})^T=-A^{(3)}_{\rm curl}$. Now it is just a matter of simple matrix multiplication 
to check that~\eqref{equation:antisymmetry} holds. Notice that for $n\geq 2$ the above reasoning can be applied row-wise.

\subsubsection{Static Maxwell equations.}\label{ex:Maxwell}
The relation between the magnetization $M:\Omega_1 \to \R^3$ of a ferromagnetic body modeled by $\Omega_1\subset\R^3$ and its 
induced magnetic field $H:\R^3\to \R^3$ 
is governed by the static Maxwell equations
\[
  \diverg(M+H) = 0,  \quad \curl H = 0\quad\text{in $\R^3$,} 
\]
(where $M$ is identified with its trivial extension by zero) or equivalently, by
\begin{align*}
 \Acal^{\rm mag}\Bigl(\begin{array}{c}M\\ H\end{array}\Bigr)=0\quad\text{in $\R^3$,}\qquad \text{with}\qquad
\Acal=\Acal^{\rm mag} =\Bigl(\begin{array}{c|c}\diverg & \diverg \\ \hline 0 & \curl\end{array}\Bigr).
\end{align*}
In fact, it is easy to verify that $\Acal^{\rm mag}$ satisfies Assumptions~\assref{ass:crp} (with $r=3$) and~\assref{ass:Ad_triangle}. 
For more details and for the precise form of $\Acal^{\rm mag}_0$ we refer to~\cite{Kre_micromag}.

Notice that problems in magnetostatics are naturally defined on the whole space.
Therefore, no extension property in the sense of Assumption~\assref{ass:approx_extension} is needed;
instead one has to adapt the reasoning of this paper to functions on $\R^3$. 
In particular, this involves replacing the Fourier series in the projection arguments of the 
next section by Fourier transforms. A detailed discussion can be found in~\cite{Kre_micromag}.
As regards Assumption~\assref{ass:antisymmetry}, it is not satisfied. 
However, when it comes to the construction of a recovery sequence in this set-up, the plane wave parts disappear anyway, 
so that Assumption~\assref{ass:antisymmetry} not needed here (see~\cite{Kre_micromag}).

\subsubsection{An example by Tartar.}\label{ex:Tartar} With $d=2$ let $u \colon \Omega_1\to \R^4$ and consider the constraint 
$\Acal u = 0$ defined through
\[
  \partial_1 u^1 + \partial_2 u^2 = 0,  \qquad \partial_1 u^3 + \partial_2 u^4 = 0,
\]
compare \cite{Tar83}. Then 
\begin{align*}
 A^{(1)}=\left(\begin{array}{cccc}
          1&0&0&0\\
          0&0&1&0
         \end{array}\right)
 \qquad\text{and}\qquad
 A^{(2)}=\left(\begin{array}{cccc}
          0&1&0&0\\
          0&0&0&1
         \end{array}\right).
\end{align*}
We point out that $\Acal$ is of constant rank with $r=2$ and trivially meets Assumption~\assref{ass:Ad_triangle}. 
Besides, one can check that Assumption~\assref{ass:antisymmetry} is satisfied as well. 
For the limit operator we have $\Acal_0 = A^{(2)}\partial_2$.
Regarding Assumption~\assref{ass:approx_extension} we observe that the trivial extension by zero of any
$u\in \ker_{\Omega_1}\Acal_0$ is $\Acal_0$-free on $\T^2$.

\subsection{Projections}
In this section an extension of the classical projection result onto $\Acal$-free vector fields of
Lemma 2.14 in~\cite{FM99} is presented.
We deal with the situation of parameter-dependent (constant-rank) operators by means of the 
Mihlin Multiplier Theorem in conjunction with a scaling argument (see Theorem~\ref{theorem_projection}). 
 
A detailed introduction to the topic of Fourier multipliers is provided for instance in \cite{Gra08} or \cite{Ste93}. 
Notice that, generally speaking, results obtained for the Fourier transform carry over to Fourier series and vice versa; for these 
transference techniques we refer to Section~3.6 of~\cite{Gra08}.

First let us recall a few basic facts from Harmonic Analysis. The Fourier transform $\Fcal f$ of a rapidly decaying test function 
$f \in \Scal(\R^d)$ is
\[
  \Fcal f(\xi):= \int_{\R^d} f(x)\ee^{-2\pi \ii x\cdot \xi}\dd{x},\qquad \xi \in \R^d,
\]
and the inverse Fourier transform is $\Fcal^{-1} f(\xi) := \Fcal f(-\xi)$. 
We also define the multiplier operator $T_m$ for $m:\R^d\setminus\{0\}\to \R$ by
\[
  T_m f=\Fcal^{-1}(m\, \Fcal f),  \qquad f\in \Scal(\R^d),
\]
whenever this makes sense. We call $m$ an $L^p$-multiplier on $\R^d$ if $T_m$ extends to a 
bounded linear operator in $L^p(\R^d)$. Finally, let $\norm{\frarg}_{\Mcal_p(\R^d)}$ denote the norm on the 
space of $L^p$-multipliers on $\R^d$,
\[
 \norm{m}_{\Mcal_p(\R^d)}=\norm{T_m}_{\Lin(L^p(\R^d);L^p(\R^d))}.
\]

For $\lambda>0$ we define a linear scaling operator $\tau_{\lambda}$ that acts on functions $m:\R^d\setminus\{0\}\to \R$ through
\[
\tau_\lambda m(\xi):=m\bigl(\xi', \lambda^{-1}\xi_d\bigr), \qquad \xi\in \R^d\setminus\{0\}.
\]
As an operator between $L^p$ spaces $\tau_\lambda:L^p(\R^d)\to L^p(\R^d)$ is bounded and its operator norm 
satisfies
\begin{align}\label{norm:tau_lambda}
\norm{\tau_\lambda}_{\Lin(L^p(\R^d);L^p(\R^d))}= \lambda^{1/p}.
\end{align}

The next lemma contains a scaling argument for Fourier multipliers on $\R^d$.
\begin{lemma}\label{lem:multiplier_scaling}
Let $\lambda>0$ and $1\leq p<\infty$. If $m:\R^d\setminus\{0\}\to \R$ is an $L^p$-multiplier, then $\tau_\lambda m$ 
is an $L^p$-multiplier as well and it holds that 
\begin{align*}
\norm{\tau_{\lambda}m}_{\Mcal_p(\R^d)} = \norm{m}_{\Mcal_p(\R^d)}.
\end{align*}
\end{lemma}

\begin{proof}
A straightforward calculation based on change of variables and the properties of the Fourier transformation 
$\Fcal$ shows that 
\[
T_{(\tau_{\lambda} m)} f = \tau_{\lambda^{-1}}T_m (\tau_\lambda f),  \qquad f \in \Scal(\R^d).
\]
Then, in view of $\eqref{norm:tau_lambda}$ one may infer 
\begin{align}\label{est:scaling}
 \norm{\tau_\lambda m}_{\Mcal_p(\R^d)} &=\norm{\tau_{\lambda^{-1}}T_m (\tau_\lambda \,\frarg)}_{\Lin(L^p(\R^d);L^p(\R^d))}\nonumber\\
&\leq \lambda^{-1/p} \;\norm{T_m}_{\Lin(L^p(\R^d);L^p(\R^d))} \;\lambda^{1/p}=\norm{m}_{\Mcal_p(\R^d)}.
\end{align}
Using the same argument once again with $\lambda$ and $m$ replaced by 
$1/\lambda$ and $\tau_\lambda m$, respectively, proves
\begin{align*}
 \norm{m}_{\Mcal_p(\R^d)} &=\norm{\tau_{1/\lambda}(\tau_\lambda m)}_{\Mcal_p(\R^d)}\leq \norm{\tau_\lambda m}_{\Mcal_p(\R^d)},
\end{align*}
which together with~\eqref{est:scaling} implies the statement.
\end{proof}

For $\eps>0$ and a constant-rank operator $\Acal$, consider for every $\xi \in \R^d \setminus \{0\}$ the orthogonal projector 
$\Pbb_\eps(\xi) \in \Lin(\R^m;\R^m)$ onto $\ker \Abb_\eps(\xi)$, and define $\Qbb_\eps(\xi) \in \Lin(\R^l;\R^m)$ by
\begin{align*}
  \Qbb_\eps(\xi)v = \begin{cases}
    z-\Pbb_\eps(\xi)z  &\text{for $v\in \range \Abb_\eps(\xi)$ with
      $v=\Abb_\eps(\xi)z$, $z\in \R^m$,} \\
    0 & \text{for $v\in \left(\range \Abb_\eps(\xi)\right)^\perp$.}
  \end{cases}
\end{align*}
Notice that $\Qbb_\eps(\xi)$ is well-defined.

Suppose that $\Acal$ satisfies Assumption~\assref{ass:crp} and let $\eps>0$.
Then, $\Pbb_\eps \colon \R^d \setminus \{0\} \to \Lin(\R^m;\R^m)$ is $0$-homogeneous and smooth, similarly 
$\Qbb_\eps \colon \R^d \setminus \{0\} \to \Lin(\R^l;\R^m)$ is $(-1)$-homogeneous and smooth. These properties 
follow from the definition and the constant-rank property of $\Acal_\eps$, which is uniform 
in $\eps$ by Lemma \ref{lemma_rank_Aeps}.

Since the operators $\Pbb_\eps$ and $\Qbb^\ast_\eps:=\Qbb_\eps(\frarg/\abs{\frarg})$ are $0$-homogeneous and smooth on 
$\Sbb^{d-1}$ for all $\eps>0$, they are
$L^p$-Fourier multipliers on $\R^d$ by the Mihlin Multiplier Theorem (see e.g.\ \cite{Gra08}) and we have the estimates 
\begin{align}
\norm{\Pbb_\eps}_{\Mcal_p(\R^d;\Lin(\R^m;\R^m))} &\leq c_d\max\{p, (p-1)^{-1}\} \, C_M(\Pbb_\eps),\label{est:Mihlin}\\
\norm{\Qbb_\eps^\ast}_{\Mcal_p(\R^d;\Lin(\R^l;\R^m))} &\leq c_d\max\{p, (p-1)^{-1}\} \, C_M(\Qbb_\eps^\ast),\label{est:Mihlin2}
\end{align}
where 
\begin{align*}
C_M(m)&= \sup\,\setb{\abs{\xi}^{\abs{\alpha}}\absb{\partial_\xi^\alpha m(\xi)}}{\xi \in \R^d\setminus \{0\}, \,
\alpha \in (\N\cup\{0\})^d \text{ such that }\abs{\alpha}\leq \lfloor d/2\rfloor +1}
\end{align*}
for any $m:\R^d\setminus\{0\}\to\R^n$ with $n\in \N$.
Notice that in \eqref{est:Mihlin} and~\eqref{est:Mihlin2} the right-hand sides still depend on $\eps$. 

Using $\Abb_\eps(\xi) = \Abb_1(\xi_\eps)$ with $\xi_\eps = (\xi',\eps^{-1}\xi_d)$, we derive the relation
\begin{align*}
 \Pbb_\eps(\xi) = \Pbb_1(\xi_\eps) = (\tau_\eps \Pbb_1)(\xi), \qquad \xi\in \R^d\setminus\{0\}.
\end{align*}
Hence, one may infer from Lemma~\ref{lem:multiplier_scaling} and $\eqref{est:Mihlin}$ that
\begin{equation} \label{eq:Pbb_eps_est}
\norm{\Pbb_\eps}_{\Mcal_p(\R^d;\Lin(\R^m;\R^m))}= \norm{\Pbb_1}_{\Mcal_p(\R^d;\Lin(\R^m;\R^m))}\leq C_{d,p,\Pbb_1}<\infty,
\end{equation}
where the upper bound is now uniform in $\eps$.
For $\Qbb_\eps^\ast$ we obtain
\begin{align*}
\Qbb_\eps^\ast(\xi) &= \Qbb_\eps\left(\frac{\xi}{\abs{\xi}}\right) 
= \frac{\abs{\xi}}{\abs{\xi_\eps}}\, \Qbb_\eps\left(\frac{\xi}{\abs{\xi_\eps}}\right)
= \frac{\abs{\xi}}{\abs{\xi_\eps}}\, \Qbb_1\left(\frac{\xi_\eps}{\abs{\xi_\eps}}\right)
= \frac{\abs{\xi}}{\abs{\xi_\eps}}\, \Qbb_1^\ast (\xi_\eps)\\
&= \frac{\abs{\xi}}{\abs{\xi_\eps}}\, (\tau_\eps \Qbb_1^\ast)(\xi) 
= m_\eps(\xi)(\tau_\eps \Qbb_1^\ast)(\xi),\qquad \xi\in \R^d\setminus\{0\}.
\end{align*}
Here $m_\eps:\R^{d}\setminus\{0\}\to \R$ is defined by $m_\eps(\xi)=\abs{\xi}/\abs{\xi_\eps}$ with 
$\xi\in \R^{d}\setminus\{0\}$. If $\eps\in (0,1]$, then
\begin{align*}
 \abs{\xi^\alpha \partial_\xi^\alpha m_\eps(\xi)} \leq 2^d
\end{align*}
for all $\xi\in (\R\setminus\{0\})^d$ and all multi-indices $\alpha\in \{0,1\}^d$. 
Hence, by the Lizorkin Multiplier Theorem \cite{Liz63}, $\{m_\eps\}_{\eps\in(0,1]}$ is a family of $L^p$-Fourier multipliers 
on $\R^d$ that is uniformly bounded in $\eps$, i.e.\
\begin{align*}
 \norm{m_\eps}_{\Mcal_p(\R^d)}\leq C_{d,p} < \infty.
\end{align*}
Together with the scaling argument of Lemma~\ref{lem:multiplier_scaling} this entails 
\begin{align}\label{est:Qbb}
\norm{\Qbb_\eps^\ast}_{\Mcal_p(\R^d;\Lin(\R^l;\R^m))}
\leq \norm{m_\eps}_{\Mcal_p(\R^d)}\norm{\tau_\eps \Qbb_1^\ast}_{\Mcal_p(\R^d;\Lin(\R^l;\R^m))}
\leq C_{d,p,\Qbb^\ast_1}<\infty
\end{align}
for all $\eps\in (0,1]$.

Finally, by transference, both $\{\Pbb_{\eps}(j)\}_{j\in \Z^d \setminus \{0\}}$ and 
$\{\Qbb_{\eps}^\ast(j)\}_{j\in \Z^d \setminus \{0\}}$ 
are discrete $L^p$-Fourier multipliers.

In the sequel, for all $\epsilon > 0$, we employ the discrete Fourier multiplier operators $\Pcal_\eps$ defined on 
$L^p(\T^d;\R^m) \cong L^p(Q^d;\R^m)$ by
\begin{equation}\label{definition_P} 
  \Pcal_\eps u(x) := \hat{u}(0) + \sum_{\xi \in \Z^d\setminus\{0\}} \Pbb_\eps(\xi) \hat{u}(\xi)
  \ee^{2\pi \ii x \cdot \xi}, \qquad x\in Q^d.
\end{equation}
Notice that we include the constant part $\hat{u}(0)$ in the definition, in contrast to other projection results.

\begin{theorem}[Projection onto $\Acal_\eps$-free fields]\label{theorem_projection} Let $p \in (1,\infty)$ and let $\Acal$ satisfy 
Assumption~\assref{ass:crp}.
Then, for every $\eps\in (0,1]$, the operators $\Pcal_\eps$ satisfy the following
properties for all $u\in L^p(\T^d;\R^m)$:
\begin{enumerate}
 \item[\it (i)] $(\Pcal_\eps \circ \Pcal_\eps) u = \Pcal_\eps u$.
 \item[\it (ii)] $\Pcal_\eps u$ is $\Acal_\eps$-free in $\T^d$.
 \item[\it (iii)] The operators $\Pcal_\eps$ are uniformly bounded with respect to $\eps$, i.e.
\[
  \qquad \norm{\Pcal_\eps u}_{L^p(\T^d;\R^m)} \leq c_p \norm{u}_{L^p(\T^d;\R^m)}
\]
with a constant $c_p > 0$ independent of $\epsilon$.
\item[\it (iv)] There exists a constant $c_p>0$ such that
\[
  \qquad \norm{u-\Pcal_\eps u}_{L^p(\T^d;\R^m)}\leq c_p \norm{\Acal_\eps u}_{W^{-1,p}(\T^d;\R^l)}
\]
for all $\eps>0$.
\item[\it (v)] Let $\eps_j\todown 0$ as $j\to \infty$ and suppose $(u_j)_j \subset L^p(\T^d;\R^m)$ is a $p$-equiintegrable sequence. 
Then the sequence $(\Pcal_{\eps_j} u_j)_j$ is still $p$-equiintegrable.
\end{enumerate}
\end{theorem}

\begin{proof} 
Assertion (iii) follows immediately from the properties of $\{\Pbb_\eps(j)\}_{j\in \Z^d\setminus \{0\}}$ 
as a discrete $L^p$-Fourier multiplier, in particular~\eqref{eq:Pbb_eps_est}.

The properties (i) and (ii) can be seen directly from the definition of $\Pcal_\eps$ in \eqref{definition_P}.
For the proof of (iv) notice that for any $u\in C^{\infty}(\T^d;\R^m)$
\[
 u(x)=\hat{u}(0) + \sum_{\xi\in \Z^d\setminus\{0\}} \hat{u}(\xi)\ee^{2\pi \ii x\cdot \xi}, \qquad x\in Q^d,
\]
where this series converges uniformly and absolutely. Then, by defining
\[
  \hat{w}_\eps(\xi):=\abs{\xi}^{-1}\Abb_\eps(\xi)\hat{u}(\xi)
\]
(notice that $\Abb_\epsilon(0)=0$ and set $\hat{w}_\eps(0)=0$)
and accounting for the $(-1)$-homogeneity of $\Qbb_\eps$, one obtains for $x\in Q^d$ that
\begin{align*}
(u-\Pcal_\eps u)(x) &= \sum_{\xi\in \Z^d\setminus\{0\}} \Qbb_\eps\bigl(\xi/\abs{\xi}\bigr)
  |\xi|^{-1} \Abb_\eps(\xi) \hat{u}(\xi)\ee^{2\pi \ii x\cdot \xi}\\
&= \sum_{\xi\in \Z^d\setminus\{0\}}\Qbb_\eps^\ast(\xi)\hat{w}_\eps(\xi)\ee^{2\pi \ii x\cdot \xi}.
\end{align*}
Since $\{\Qbb_\eps^\ast(j)\}_{j\in \Z^d\setminus\{0\}}$ are discrete $L^p$-Fourier multipliers with norms bounded uniformly 
with respect to $\eps\in (0,1]$, see~\eqref{est:Qbb}, we have
\[
\|u-\Pcal_\eps u\|_{L^p(\T^d;\R^m)}\leq c_p \|w_\eps\|_{L^p(\T^d;\R^l)}.
\]
Using the definitions of Sobolev spaces on the torus, see Section~\ref{sec:Sobolev_torus}, and the fact 
that $\limsup_{\abs{\xi}\to \infty} (1+4\pi^2\abs{\xi}^2)^{1/2}/\abs{\xi} < \infty$ we have (again by the Mihlin Multiplier Theorem) that
\[
  \norm{w_\eps}_{L^p(\T^d;\R^l)} \leq c_p \norm{\Acal_\eps u}_{\Wrm^{-1,p}(\T^d;\R^l)},
\]
and so the claim holds in the case of smooth functions. The general result for $u\in L^p(\T^d;\R^m)$ follows by a density argument. 

In view of (iii) the proof of (v) is exactly the same as the one of Lemma~2.14~(iv) in~\cite{FM99}. 
\end{proof}

\begin{remark} 
The essential improvement of Theorem~\ref{theorem_projection} in comparison to Lemma~2.14 in \cite{FM99} is that all 
constants are uniform with respect to $\eps$. 
In the above reasoning we employed a scaling argument together with the Mihlin Multiplier Theorem, before arguing with the Lizorkin Multiplier Theorem, to achieve that.
Alternatively, one can use the Lizorkin Multiplier Theorem~\cite{Liz63} right away.
Then the expressions $C_M(\Pbb_\eps)$ in $\eqref{est:Mihlin}$ and $C_M(\Qbb_\eps^\ast)$ in $\eqref{est:Mihlin2}$ are replaced by
\begin{align*}
C_L(\Pbb_\eps)&= \sup\,\setb{\absb{\xi^\alpha\partial_\xi^\alpha \Pbb_\eps(\xi)}}{\xi \in (\R\setminus \{0\})^d, \,
\alpha \in \{0,1\}^d}\\
&=\sup\,\setb{\absb{\xi^\alpha\partial_\xi^\alpha \Pbb_1(\xi)}}{\xi \in (\R\setminus \{0\})^d, \,
\alpha \in \{0,1\}^d} = C_L(\Pbb_1),
\end{align*}
and
\begin{align*}
 C_L(\Qbb^\ast_\eps)&= \sup\,\setb{\absb{\xi^\alpha\partial_\xi^\alpha \Qbb^\ast_\eps(\xi)}}{\xi \in (\R\setminus \{0\})^d, \,
 \alpha \in \{0,1\}^d}\\
& \leq\,\sup\,\setb{\absb{\xi^\alpha\partial_\xi^\alpha m_\eps(\xi)}}{\xi \in (\R\setminus \{0\})^d, \,
 \alpha \in \{0,1\}^d}\\ 
& \qquad\cdot\sup\,\setb{\absb{\xi^\alpha\partial_\xi^\alpha \Qbb^\ast_1(\xi)}}{\xi \in (\R\setminus \{0\})^d, \,
 \alpha \in \{0,1\}^d}\\
 & = C_L(m_\eps) \, C_L(\Qbb_1^\ast) \leq 2^d\, C_L(\Qbb_1^\ast),
\end{align*}
respectively, for $\eps\in(0,1]$.
\end{remark}

\subsection{$\Acal_\eps$-quasiconvexity and asymptotic $\Acal_0$-quasiconvexity}\label{sec:A-quasiconvex}

The notion of $\Acal$-quasicon\-vexity was first introduced and studied by Dacorogna \cite{Dac82}.

\begin{definition}[$\Acal$-quasiconvexity]\label{def:Aqc}
A function $f\colon \R^m \to \R$ is called $\Acal$-quasiconvex, if
\begin{align*}
 f(v)\leq \dashint_{Q^d} f(v+w(y))\dd{y}
\end{align*}
for all $v\in \R^m$ and all $w\in C^{\infty}(\T^d;\R^m)$ with $\Acal w=0$ in $\T^d$ and $\int_{Q^d} w \dd y = 0$.
\end{definition}

\begin{remark}\label{remark_A-quasiconvexity}
1. By a simple scaling argument it can be seen that the choice of $Q^d$ as a domain is not essential, 
but can be replaced by any open cuboid $E\subset\R^d$, if we choose test functions 
$w\in C^{\infty}(\T^d(E);\R^m)$ with $\Acal w=0$ in $\T^d(E)$ and $\int_{E} w \dd y = 0$.

2. If $f$ is a continuous function with $p$-growth, i.e.\ $\abs{f(v)}\leq C(1+\abs{v}^p)$ for all $v\in \R^m$,
the space $C^\infty(\T^d;\R^m)$ in Definition~\ref{def:Aqc} may be replaced by $L^p(Q^d;\R^m)$ (compare Remark~3.3~(ii) of \cite{FM99}).

3. Let us point out that the above definition does not need the operator $\Acal$ to satisfy a constant-rank property.
\end{remark}

The $\Acal$-quasiconvex envelope of a function $f\colon\R^m \to \R$, denoted $\Qcal_\Acal f$,
is defined as
\[
\Qcal_\Acal f(v) = \inf\set{\dashint_{Q^d} f(v+ w(y))\dd{y}}{ w\in C^{\infty}(\T^d;\R^m)\cap \ker_{\T^d}\Acal,  
\textstyle\,\int_{Q^d} w \dd y = 0}.
\]
If $\Acal$ is of constant rank 
and $f$ is continuous, then $\Qcal_\Acal f$ can be proven to be $\Acal$-quasiconvex and 
upper semicontinuous (see Proposition~3.4 of \cite{FM99}).
Notice that 
in general, though, $\Qcal_\Acal f$ is not continuous even if $f$ is smooth. Counterexamples can be found in 
Remark~3.5\,(ii) of \cite{FM99}.
In the special cases $\Acal=\diverg$ and $\Acal=\curl$, however, continuity ensues. Indeed, $\Qcal_{\diverg} f$ is exactly the 
convexification of $f$, while $\Qcal_{\curl} f$ is a quasiconvex function and hence continuous.

It is instructive to observe that the notion of $\Acal_\eps$-quasiconvexity is independent of $\eps$:

\begin{lemma}\label{lem:Aeps-Adelta_quasiconvexity}
For all $\eps>0$ and $\delta>0$ it holds that a function $f\colon\R^m\to \R$ is $\Acal_\eps$-quasiconvex 
if and only if it is $\Acal_\delta$-quasiconvex.
\end{lemma}

As a consequence of the lemma, $\Qcal_{\Acal_\eps}f=\Qcal_{\Acal_\delta}f = \Qcal_\Acal f$ for all $\eps, \delta>0$ (notice $\Acal_1 = \Acal$).

\begin{proof}
Assume that $f$ is $\Acal_\eps$-quasiconvex and let $w\in C^{\infty}(\T^d;\R^m)$ with $\Acal_\delta w=0$ in 
$\T^d$, $\int_{Q^d} w \dd x = 0$. 
Consider the transformation of variables given by $z=(z',z_d)=(y', (\delta/\eps) y_d)$ and set 
$\tilde{w}(z)=w(z', (\eps/\delta) z_d)$ for $z\in E := (0,1)^{d-1} \times (0, \delta / \eps)$. 
Then, $\tilde{w}$ is $E$-periodic, $\Acal_\eps$-free in $\T^d(E)$ and has mean value zero over $E$. 
Hence, in view of Remark~\ref{remark_A-quasiconvexity}~1.~it follows by the $\Acal_\eps$-quasiconvexity of $f$ that 
\begin{align*}
\dashint_{Q^d} f(v+w(y))\dd{y}
&=\frac{\eps}{\delta}\int_{E} f\big(v+w\left(z', (\eps/\delta) z_d\right)\big)\dd{z} \\
&= \dashint_{E} f(v+\tilde{w}(z))\dd{z} \geq f(v)
\end{align*}
for all $v\in \R^m$. 
\end{proof}

As already introduced in the introduction, we also define the asymptotic $\Acal_0$-quasiconvex envelope of a continuous $f:\R^m\to \R$ as follows:
\[
  \Qcal^\infty_{\Acal_0} f(v) = \lim_{\eta \to \infty} \Qcal^\eta_{\Acal_0} f(v) = \sup_{\eta > 0} \Qcal^\eta_{\Acal_0} f(v), \qquad v\in \R^m,
\]
where
\begin{align}\label{def:Qcaleta}
\Qcal^\eta_{\Acal_0} f(v) &:= \inf\biggl\{ \dashint_{Q^d} f(v+ w(y))\dd{y}\ \textup{\textbf{:}}\ 
 w \in C^{\infty}(\T^d;\R^m),
\\ &\qquad\qquad\qquad\qquad\qquad \eta \norm{\Acal_0 w}_{W^{-1,1}(\T^d;\R^l)}\leq 1,\;\; \textstyle\int_{Q^d} w \dd y = 0 \biggr\}\nonumber
\end{align}
with $W^{-1,1}(\T^d;\R^l)$ the dual space of $W_0^{1,\infty}(\T^d;\R^l)$.

Analogously to Proposition~3.4 in~\cite{FM99}, we get that $\Qcal^\eta_{\Acal_0} f$ is upper semicontinuous for all $\eta > 0$. 
Moreover, if $f$ is continuous with $p$-growth, an approximation argument similar to that in Remark~\ref{remark_A-quasiconvexity} 2. allows us to replace $C^{\infty}(\T^d;\R^m)$ in \eqref{def:Qcaleta} with $L^p(Q^d;\R^m)$. Also notice that by definition $\Qcal^\infty_{\Acal_0}f\leq \Qcal_{\Acal_0}f$.

\begin{remark} \label{rem:envelopes}
If $\Acal$ and $f:\R^m\to \R$ meet the requirements of Theorem~\ref{mainresult_Gamma}, one obtains that
$\int_{\Omega_1}\Qcal^\infty_{\Acal_0} f(u)\dd{x}\leq \int_{\Omega_1}\Qcal_\Acal f(u) \dd{x}$ 
for all $u\in L^p(\Omega_1;\R^m)\cap \ker_{\Omega_1}\Acal_0$. Testing with constant fields results in
$\Qcal^\infty_{\Acal_0} f(v) \leq \Qcal_\Acal f(v)$ for all $v\in \R^m$. Hence, if $f$ is asymptotically $\Acal_0$-quasiconvex, i.e.~$f=\Qcal^\infty_{\Acal_0} f$,
we may argue that $f=\Qcal^\infty_{\Acal_0}f\leq \Qcal_\Acal f\leq f$, 
which implies $\Qcal_\Acal f=f$ and 
in particular the $\Acal$-quasiconvexity of $f$.
\end{remark}

We say that a function $f:\Omega_1\times\R^m\to \R$ 
is $\Acal$-quasiconvex, 
if $f(x,\frarg)$ is $\Acal$-quasiconvex for almost all $x\in \Omega_1$. Accordingly, the $\Acal$-quasiconvex and asymptotic $\Acal_0$-quasiconvex envelopes of $f$, 
that is $\Qcal_\Acal f, \Qcal_{\Acal_0}^\infty f:\Omega_1\times \R^m\to \R$, are given by $\Qcal_\Acal f(x,\frarg)$ and $\Qcal_{\Acal_0}^\infty f(x, \frarg)$, respectively, in the above sense for almost every $x\in \Omega_1$.

\begin{lemma}\label{lem:normal_integrand}
Let $\Acal$ have constant rank.
If $f:\Omega_1\times\R^m\to \R$ is a Carath\'eodory function, then $-\Qcal_\Acal f$ is a normal integrand.
\end{lemma}

\begin{proof}
In view of the upper semicontinuity of $\Qcal_\Acal f(x, \frarg)$ for almost all $x\in \Omega_1$ and of 
Theorem~6.28 of~\cite{FL07} it is sufficient to show that for every $\epsilon>0$ there 
exists a closed set 
$K_\epsilon\subset \Omega_1$ with $\abs{\Omega_1\setminus K_\epsilon}\leq \epsilon$ such that $\Qcal_\Acal f|_{K_\epsilon\times\R^m}$ 
is upper semicontinuous.

For $K_\epsilon$ we pick the compact set resulting from the Scorza--Dragoni Theorem applied to $f$ and proceed similarly to 
the proof of Proposition~3.4 (Case 1) in \cite{FM99} 
showing that $\Qcal_\Acal^R f$, which is defined for almost every $x\in \Omega_1$ and 
$v\in \R^m$ by
\begin{align*}
\Qcal_\Acal^R f(x,v) :=\inf \bigg\{ \, \dashint_{Q^d}f(x,v+w(y))\dd{y} \ \textup{\textbf{:}}\  
&w\in C^{\infty}(\T^d;\R^m) \cap \ker_{\T^d} \Acal, \\ 
&\textstyle\int_{Q^d}w\dd{y} = 0, \, \norm{w}_{L^{\infty}(Q^d;\R^m)}\leq R \, \bigg\} 
\end{align*}
is continuous on $K_\epsilon\times \R^m$ for all $R>0$. If $(x_k)_k \subset K_\epsilon$ with $x_k\to x$ and 
$(v_k)_k\subset \R^m$ with $v_k\to v$, then 
\begin{align*}
\limsup_{k\to \infty} \Qcal_\Acal f(x_k, v_k)\leq \limsup_{k\to \infty} \Qcal_\Acal^R f(x_k, v_k)= \Qcal_\Acal^R f(x, v).
\end{align*}
Finally, since $(\Qcal_\Acal^R f)_R$ is a decreasing sequence converging pointwise to $\Qcal_\Acal f$ as 
$R\to \infty$, the claim is proven.
\end{proof}

\begin{remark}\label{rem:measurability} 
1. The preceding lemma implies that for a constant-rank operator $\Acal$
the integrand $\Qcal_\Acal f(\frarg, u(\frarg))$ 
in Theorem~\ref{mainresult_Gamma} is Lebesgue measurable.
We remark that measurability of $\Qcal_\Acal f(\frarg, u(\frarg))$, where $u$ is an $\Acal$-free $L^p$-function in $\Omega_1$ and 
$f$ satisfies $\eqref{f_growth/coercivity}$, also follows implicitly from the localization approach in the proof of the 
relaxation result Theorem~1.1 in~\cite{BFL00} and is essentially a consequence of the representation formula in Lemma~3.5 
of~\cite{BFL00} in conjuction with the Radon--Nikod\'ym Theorem.

2. Notice that if $f$ is Carath\'{e}odory, then for $u \colon \Omega_1 \to \R^m$ the compound function $\Qcal^\infty_{\Acal_0} f(\frarg, u(\frarg))$ is measurable since it is the pointwise limit of the measurable functions $\Qcal^\eta_{\Acal_0} f(\frarg, u(\frarg))$ for a sequence $\eta \to \infty$ (by an argument analogous to Lemma~\ref{lem:normal_integrand}, $-\Qcal^\eta_{\Acal_0} f$ is a normal integrand).
\end{remark}


\section{Proof of the lower bound}\label{sec:proof_lower-bound}

To show the lower bound we follow a classical Young measure approach, which requires two important 
technical tools, a decomposition lemma 
and a localization result. Let us remark that the proofs in this section are all of local nature and 
therefore do not require the existence 
of approximate extensions as in Assumption~\assref{ass:approx_extension}. A major difficulty results from the fact that owing to the non-constant rank nature of general $\Acal_0$, a projection operator onto $\Acal_0$-free vector fields satisfying good estimates (in the sense of Theorem~\ref{theorem_projection}) does not seem to exist in general.


\subsection{Equiintegrability}

The decomposition lemma of Fonseca and M{\"u}ller~\cite{FM99} generalizes the original work in the 
gradient setting (cf.\ Lemma~1.2 of~\cite{FMP98}, a similar result was also obtained by Kristensen~\cite{Kri99}) to 
the context of $\Acal$-free fields.

\begin{lemma}[adapted from~{Lemma 2.15 of~\cite{FM99}}]\label{lem:decomposition_FM99} 
Let $1\leq q<p < \infty$ and suppose that $(u_j)_j$ is a bounded sequence in $L^p(Q^d;\R^m)$ with 
$u_j\weakly u$ in $L^p(Q^d;\R^m)$ and $\Acal u_j \to 0$ in $W^{-1,p}(Q^d;\R^l)$. Then, $u_j$ can be decomposed as
\begin{align*}
u_j=z_j+r_j, 
\end{align*}
where $(z_j)_j \subset L^p(Q^d; \R^m)$ is a $p$-equiintegrable sequence that satisfies
\[
  \Acal z_j=0  \quad\text{in $\T^d$, }\qquad
  \int_{Q^d} z_j\dd{y}= \int_{Q^d} u\dd{y}, \qquad\text{for all $j \in \N$,}
\]
and $(r_j)_j \subset L^p(Q^d;\R^m)$ is such that $r_j\to 0$ in $L^q(Q^d;\R^m)$ as $j\to \infty$.
\end{lemma}

Equiintegrability in the context of thin films with functionals depending on gradients was first studied in~\cite{BF02},
an elegant alternative proof is given in \cite{BZ07}. 
The following result is the appropriate generalization of Lemma~\ref{lem:decomposition_FM99} to the context of dimension reduction 
problems with parameter-dependent operators. It involves the essential cut-off that is needed for the localization procedure 
in Proposition~\ref{prop_localization}. 

\begin{theorem}[Decomposition lemma]\label{prop:decomposition2}
Let $1\leq q<p<\infty$ and ${\eps_j}\todown 0$ as $j \to \infty$. Further, suppose that $\Acal$ 
satisfies Assumptions~\assref{ass:crp} and \assref{ass:Ad_triangle}, and assume 
that $(u_j)_j$ is a bounded sequence in $L^p(Q^d;\R^m)$ with 
\begin{align*}
  \text{$u_j \weakly 0$ in $L^p(Q^d;\R^m)$} \quad\text{and}\quad
  \text{$\Acal_{\eps_j} u_j \to 0$ in $W^{-1,p}(Q^d;\R^l)$.}
\end{align*}
Then,
\begin{align*}
 u_j=w_j +r_j,
\end{align*}
where $(w_j)_j \subset L^p(Q^d; \R^m)$ is a $p$-equiintegrable sequence that satisfies 
\[
\Acal_0 w_j\to0 \text{ in $W^{-1,q}(\T^d;\R^l)$}, \qquad \int_{Q^d} w_j\dd{y}= 0
\text{ for all $j \in \N$,}
\]
and $(r_j)_j\subset L^p(Q^d;\R^m)$ is such that $r_j\to 0$ in $L^q(Q^d;\R^m)$ as $j\to \infty$. 
\end{theorem}

%
%

\begin{proof}
We observe that $\Acal_{\eps_j} u_j \to 0$ in $W^{-1,p}(Q^d;\R^l)$ implies
\begin{align}\label{convergenceA0uj}
 \Acal_0 u_j\to 0 \qquad \text{in $W^{-1,p}(Q^d;\R^l)$}
\end{align}
as $j\to \infty$. Indeed, with the notation of Section~\ref{sec:operator_A}, in particular~\eqref{notationA+}, and
\begin{align*}
 \Acal_0=\left[\begin{array}{c}(\Acal_0)_+\\ \hline (\Acal_0)_-\end{array}\right],
\end{align*}
it holds that
\begin{align*}
 (\Acal_0)_+ u_j = A^{(d)}_+\partial_d u_j= \eps_j (\Acal_{\eps_j})_+ u_j - \eps_j \Acal'_+ u_j\ \to\  0 
\quad\text{ in $W^{-1,p}(Q^d;\R^r)$,}
\end{align*}
since $\norm{\Acal'_+ u_j}_{W^{-1,p}(Q^d;\R^r)}\leq c\,\norm{u_j}_{L^p(Q^d;\R^m)}\leq c$, and
\begin{align*}
 (\Acal_0)_- u_j = \Acal'_- u_j =(\Acal_{\eps_j})_- u_j \to\  0\ \qquad\text{ in $W^{-1,p}(Q^d;\R^{l-r})$}
\end{align*}
as $j\to \infty$.
This proves \eqref{convergenceA0uj}.
By truncation one finds a $p$-equiintegrable sequence $(z_j)_j$ with $z_j - u_j\to 0$ in $L^q(Q^d;\R^m)$.
Then, 
\begin{align}\label{conv_A0}\Acal_0 z_j\to 0 \qquad \text{in $W^{-1,q}(Q^d;\R^l)$,}\end{align} and $z_j\to 0$ in $W^{-1, p}(Q^d;\R^m)$ as $j\to \infty$. 
The latter results from $z_j\weakly 0$ in $L^p(Q^d;\R^m)$ together with the compact embedding 
$L^p(Q^d;\R^m)\hookrightarrow W^{-1, p}(Q^d;\R^m)$.

Performing a suitable cut-off allows us to switch from $W^{-1,q}(Q^d;\R^l)$ to $W^{-1,q}(\T^d;\R^l)$ in \eqref{conv_A0}. 
For all $\varphi\in C_c^\infty(Q^d;[0,1])$ we argue that
\begin{align}\label{eq:cut-off}
 \Acal_{0}(\varphi z_j) = \varphi \Acal_{0}z_j + \sum_{k=1}^{d} A_0^{(k)} z_j \partial_k \varphi \ \to\  0 
\quad\mbox{in $W^{-1,q}(\T^d;\R^l)$ as $j\to \infty$}.
\end{align}
Indeed, the convergence of the first term is due to the sequence of truncated functions $(z_j)_j$ satisfying 
\eqref{conv_A0} in combination with the fact that $\varphi$ has compact support in $Q^d$. 
Moreover, since $z_j\to 0$ in $W^{-1, p}(Q^d;\R^m)$, also the second term in~\eqref{eq:cut-off} converges to zero.

We may now pick a sequence $(\varphi_j)_j\subset C_c^\infty(Q^d;[0,1])$ of cut-off functions with
$\varphi_j\to 1$ for $j\to \infty$ such that $(v_j)_j\subset L^p(Q^d;\R^m)$ defined by $v_j=\varphi_j z_j$ is a 
$p$-equiintegrable sequence and satisfies
\begin{align*}
 v_j\weakly 0\quad\text{in $L^q(Q^d;\R^m)$} \qquad \text{and}\qquad \Acal_0 v_j \to 0 \quad\text{in $W^{-1,q}(\T^d;\R^l)$}.
\end{align*}
Finally, setting 
$w_j=v_j -\int_{Q^d}v_j\dd{y}$ and $r_j=u_j-w_j$ for all $j\in \N$ 
provides the sought sequences.
\end{proof}


\subsection{Localization via Young measures}
We start this section by presenting a formulation of the fundamental theorem on Young measures that will be needed in the sequel. 
For more general statements and proofs the reader is referred for instance to~\cite{Mue99, Ped97}.

\begin{theorem}[Fundamental Theorem on Young measures]\label{theorem_YM}
Suppose $U \subset \R^d$ is an open, bounded set and let the sequence $(z_j)_j$ be bounded in $L^1(U;\R^m)$. Then there exist a 
subsequence $(z_j)_j$ (not relabeled) and a weak$\ast$ measurable map $\nu: U \to \Mcal_1(\R^m)$, where $\Mcal_1(\R^m)$ denotes 
the space of probability measures on $\R^m$, such that the following holds:
\begin{enumerate}
 \item[\it (i)] For all $g \in C_0(\R^m)$
\[
  \qquad\qquad g(z_j) \weaklystar \bigl(x\mapsto \dpr{g,\nu_x}\bigr)\text{ in } L^{\infty}(U),
  \qquad\text{where}\qquad
  \dpr{g,\nu_x}=\int_{\R^m}g(y)\dd{\nu_x(y)}. 
\]
\item[\it (ii)] Let $f\in C(\R^m)$. Then 
\[
 f(z_j) \weakly \bigl(x\mapsto \dpr{f,\nu_x}\bigr)\text{ in } L^{1}(U)
 \qquad\text{if $\left(f(z_j)\right)_j$ is equiintegrable.}
\]
\item[\it (iii)] If $f:U\times \R^m\to \R$ is Carath{\'e}odory and bounded from below, then 
\[
\liminf_{j\to \infty} \int_{U} f\bigl(x,z_j(x)\bigr)\dd{x}\geq \int_U \dpr{f(x,\frarg),\nu_x} \dd{x}.
\]
\end{enumerate}
\end{theorem}

The map $\nu$ is called the Young measure generated by the sequence $(z_j)_j$. We will be using the shorthand notation
\[
 z_j \toYM (\nu_x)_{x\in U}.
\]

A proof of the following result can be found in Proposition 2.4 of~\cite{FM99}.

\begin{lemma}\label{lemma_YM}
With $U$ as in Theorem~\ref{theorem_YM}, let $(z_j)_j$ and $(w_j)_j$ be bounded sequences in $L^1(U;\R^m)$ such that 
$(z_j)_j$ generates the Young measure $(\nu_x)_{x\in U}$ and $w_j\to 0$ in measure for $j\to \infty$. 
Then
\begin{align*}
z_j+w_j\toYM (\nu_x)_{x\in U}.
\end{align*}
\end{lemma} 

Next, we employ the blow-up technique to prove the following localization result, which is necessary for 
obtaining Jensen-type inequalities. 
These in turn will then imply the liminf-inequality. Notice that in comparison to Proposition 3.8 in \cite{FM99}, we do not 
need $p$-equi\-integrability of $(u_j)_j$.

\begin{proposition}[Localization]\label{prop_localization}
Let 
$1\leq q < p < \infty$. 
Suppose $(u_j)_j\subset L^p(\Omega_1;\R^m)$ is a sequence such that $u_j\weakly u$ in $L^p(\Omega_1;\R^m)$ and 
$\Acal_{\eps_j} u_j \to 0$ in $W^{-1,p}(\Omega_1;\R^l)$ as $j\to \infty$, where $\eps_j\todown 0$. 
Further, let $(\nu_x)_{x\in \Omega_1}$ be the Young measure generated by $(u_j)_j$. 
Then, for almost every $a\in \Omega_1$ there exist a subsequence of $(\eps_j)_j$ (not relabeled) and a $p$-equiintegrable sequence 
$(z_j)_j\subset L^p(Q^d;\R^m)$ with $\Acal_0 z_j\to 0$ in $W^{-1,q}(\T^d;\R^l)$ such that
\[
\int_{Q^d}z_j\dd{y}=u(a)\qquad \text{for $j\in \N$,}
\] 
and $(z_j)_j$ generates the homogeneous Young measure $(\nu_a)_{y\in Q^d}$, i.e.\ $z_j \toYM (\nu_a)_{y\in Q^d}$.
\end{proposition}

\begin{proof}
Assume $\Lfrak$ and $\Cfrak$ are countable dense subsets of 
$L^1(Q^d)$ and $C_0(\R^m)$, respectively, which determine the Young measure convergence, i.e.\ 
for a sequence $(v_j)_j \subset L^p(Q^d;\R^m)$ the validity of
\[
  \int_{Q^d} \psi(y) g(v_j(y)) \dd{y} \to \int_{Q^d} \psi(y) \, \dpr{g, \mu_y} \dd{y}
  \qquad\text{for all $\psi \in \Lfrak$, $g \in \Cfrak$}
\]
with a weak$\ast$ measurable $\mu:Q^d\to \Mcal_1(\R^m)$ implies $v_j \toYM (\mu_y)_{y \in Q^d}$. 
Without loss of generality the elements of $\Lfrak$ are smooth up to the boundary of $Q^d$. 

By $\Omega_0 \subset \Omega_1$ we denote the set of points $a\in \Omega_1$ that are Lebesgue points for 
the functions $x\mapsto \dpr{\id,\nu_x}=u(x)$, $x\mapsto \dpr{g,\nu_x}$, and $x\mapsto\int_{\R^d}\abs{\frarg}^p \dd \nu_x$. In particular,
\begin{equation}\label{Lebesgue_point2}
 \lim_{R\to 0} \int_{Q^d} \bigl| \dpr{g,\nu_{a+Ry}}-\dpr{g,\nu_a}\bigr| \dd{y}=0
 \qquad\text{for all $g \in \Cfrak$.}
\end{equation}
By standard results in measure theory, $\abs{\Omega_0} = \abs{\Omega_1}$.

Now fix $a\in \Omega_0$. For $R>0$ sufficiently small, i.e.\ $R$ small enough such that $a+RQ^d \subset\subset\Omega_1$, we define 
\[
 u_{R,j}(y)= u_j(a+Ry) -\dpr{\id,\nu_a},\ \ \ y\in Q^d.
\]
Then, Theorem \ref{theorem_YM}~(i) in combination with Lebesgue's convergence theorem and $\eqref{Lebesgue_point2}$ implies 
for all $\psi\in \Lfrak$ and all $g\in \Cfrak$ that
\begin{align}\label{limit_w1}
&\lim_{R\to 0} \lim_{j\to \infty} \int_{Q^d} \psi(y) g\bigl(u_{R,j}(y)+\dpr{\id,\nu_a}\bigr) \dd{y} \nonumber
=  \lim_{R\to 0} \lim_{j\to \infty} \int_{Q^d} \psi(y) g\bigl(u_j(a+Ry)\bigr) \dd{y}\\
&\qquad=\lim_{R\to 0} \int_{Q^d} \psi(y) \dpr{g,\nu_{a+Ry}} \dd{y} =\dpr{g,\nu_a} \int_{Q^d} \psi(y)\dd{y}.
\end{align}
By a similar argument we deduce that
\[
\lim_{R\to 0}\lim_{j\to \infty} \norm{u_{R,j}+\dpr{\id,\nu_a}}_{L^p(Q^d;\R^m)}
\]
exists and is finite. Observe also that $u_{R,j}\weakly 0$ in $L^p(Q^d;\R^m)$ as $j\to \infty$ and $R\to 0$ (in this order), by an argument analogous 
to~\eqref{limit_w1} with $g = \id$.
Next, we will show that
\begin{equation}\label{convergence_Aepsw}
\lim_{R\to 0}\lim_{j\to \infty} \|\Acal_{\eps_j}u_{R,j}\|_{W^{-1,p}(Q^d;\R^l)}=0.
\end{equation}
If not stated otherwise, the supremum in the following estimate is taken with respect to $v\in W_0^{1,p'}(Q^d;\R^l)$ with 
$\|v\|_{W^{1,p'}(Q^d;\R^l)}\leq 1$. So,
\begin{align*}
&\norm{\Acal_{\eps_j}u_{R,j}}_{W^{-1,p}(Q^d;\R^l)}  = \sup\, \absBB{\int_{Q^d}  u_j(a+Ry)  \cdot \Acal_{\eps_j}^T v(y)\dd{y}}\nonumber\\
& = R^{-d} \, \sup\, \absBB{\int_{a+RQ^d}  u_j(x)  \cdot (\Acal_{\eps_j}^T v) \Bigl(\frac{x-a}{R}\Bigr)\dd{x}}\nonumber\\
& \leq c_R \, \sup\,\set{  \absBB{\int_{\Omega_1}  u_j  \cdot \Acal_{\eps_j}^T z\dd{x}}}
{z\in W_0^{1,p'}(\Omega_1;\R^l),\|z\|_{W^{1,p'}(\Omega_1;\R^l)}\leq 1}\nonumber\\
&=  c_R \, \norm{\Acal_{\eps_j} u_j}_{W^{-1,p}(\Omega_1;\R^l)}.
\end{align*}
Taking the limit $j\to\infty$ makes the above expression tend to zero for any fixed $R$ and $\eqref{convergence_Aepsw}$ is proven.

In view of $\eqref{limit_w1}$ and $\eqref{convergence_Aepsw}$ we can finally extract a diagonal sequence 
$(u_k)_{k}\subset L^p(Q^d;\R^m)$ that satisfies
\begin{enumerate}
\item[(i)]   $u_k\weakly 0$ in $L^p(Q^d; \R^m)$,
\item[(ii)]  $\Acal_{\eps_k} u_k \to 0$ in $W^{-1,p}(Q^d; \R^l)$,
\item[(iii)] $\int_{Q^d} \psi(y) g(u_k(y)+\dpr{\id,\nu_a}) \dd{y} \to \dpr{g, \nu_a}\int_{Q^d} \psi(y) \dd{y}$
  for all $\psi \in \Lfrak$, $g \in \Cfrak$.
\end{enumerate}

By Theorem~\ref{prop:decomposition2}
there exists a $p$-equiintegrable sequence $(z_k)_{k}\subset L^p(Q^d;\R^m)$ with $\Acal_0 z_k \to 0$ in $W^{-1,q}(\T^d;\R^l)$ and 
\begin{align}\label{convergence_wk}
 \|u_k+\dpr{\id,\nu_a}-z_k\|_{L^q(Q^d;\R^m)}\to 0 \qquad\text{ as $j\to \infty$},
\end{align}
satisfying $\int_{Q^d} z_k \dd{y}= \dpr{\id,\nu_a} = u(a)$.  Finally, (iii) implies that the homogeneous Young measure 
$(\nu_a)_{y\in Q^d}$ is generated by the 
sequence $(u_k+\dpr{\id,\nu_a})_{k}$, so that in view of $\eqref{convergence_wk}$ and Lemma~\ref{lemma_YM},
\begin{equation*}
z_k\toYM (\nu_a)_{y\in Q^d}.
\end{equation*}
This concludes the proof.
\end{proof}

\subsection{Liminf-inequality}

Finally, we are in the position to prove the $\liminf$-inequality of Theorem~\ref{mainresult_Gamma}.


\begin{proof}[Proof of Theorem~\ref{mainresult_Gamma}~(i)]
The inclusion $u\in \Ucal_0$ is easy to show,
just compare the definitions of $\Acal_\eps$ and $\Acal_0$ in $\eqref{Aeps_A0}$, and use the weak $L^p$-convergence of $u_j$ to $u$.
After extracting a subsequence (not relabeled) we may further assume that 
$\liminf_{j\to \infty} F_{\eps_j}[u_j] = \lim_{j\to\infty} F_{\eps_j}[u_j]<\infty$ and $u_j\toYM (\nu_x)_{x\in \Omega_1}$.

Since by assumption $f$ is Carath{\'e}odory and bounded from below in view of \eqref{f_growth/coercivity}, Theorem~\ref{theorem_YM}~(iii) yields
\begin{equation}\label{limit3}
\liminf_{j\to \infty}  \int_{\Omega_1} f\bigl(x,u_j(x)\bigr)\dd{x}  
\geq \int_{\Omega_1} \dpr{f(x,\frarg),\nu_x} \dd{x}.
\end{equation}
The essential step is now to derive the appropriate Jensen-type inequalities for almost every $a \in \Omega_1$. 
As we will see, they follow from the localization principle proved in Proposition~\ref{prop_localization} in 
conjunction with the properties of the asymptotic $\Acal_0$-quasiconvex envelope of $f$.

From Proposition~\ref{prop_localization} we obtain for almost every $a\in \Omega_1$ a $p$-equiintegrable sequence 
$(z_j)_j\subset L^p(Q^d;\R^m)$ with $\Acal_0 z_j\to 0$ in $W^{-1,1}(\T^d;\R^l)$ that generates the homogeneous Young measure $(\nu_a)_{y\in Q^d}$ 
and satisfies $\int_{Q^d}z_j\dd{y}=u(a)$. Let us fix such an $a\in \Omega_1$ with $f(a,\frarg)\in C(\R^m)$.
By the growth conditions on $f$, the sequence $\bigl( f\bigl(a, z_j\bigr)\bigr)_j$ is equiintegrable, so that the 
fundamental theorem on Young measures, Theorem~\ref{theorem_YM}~(ii), gives
\begin{equation}\label{limit1}
\dpr{f(a,\frarg),\nu_a}=\lim_{j\to \infty} \int_{Q^d} f\bigl(a, z_j(y)\bigr)\dd{y}.
\end{equation}
Accounting for the properties of ${\Qcal^\infty_{\Acal_0}} f(a,\frarg)$ 
leads to
\begin{align}
\liminf_{j\to \infty} \int_{Q^d}  f\bigl(a,z_j(y)\bigr)\dd{y} 
&=\liminf_{j\to \infty} \dashint_{Q^d} f\bigl(a,z_j(y)-u(a)+u(a)\bigr)\dd{y} \nonumber\\
&\geq \lim_{\eta\to\infty} \Qcal^\eta_{\Acal_0} f\bigl(a,u(a)\bigr) \nonumber\\
&\geq\Qcal^\infty_{\Acal_0} f\bigl(a,u(a)\bigr).  \label{limit2}
\end{align}
Combining $\eqref{limit1}$ and $\eqref{limit2}$, we find
\[
 \dpr{f(a,\frarg),\nu_a} \geq \Qcal^\infty_{\Acal_0} f\bigl(a,u(a)\bigr).
\]
Consequently, in view of $\eqref{limit3}$ this implies
\begin{equation*}
\liminf_{j\to \infty} \int_{\Omega_1}  f\bigl(x,u_j(x)\bigr)\dd{x} \geq \int_{\Omega_1} \Qcal^\infty_{\Acal_0} f\bigl(x,u(x)\bigr)\dd{x},
\end{equation*}
and the assertion follows.
\end{proof}


\section{Proof of the upper bound}\label{sec:proof_recovery}

For the proof of the upper bound we proceed in two steps. First, by using the relaxation result of Theorem~1.1 in~\cite{BFL00} 
we may reduce our considerations to a functional whose integrand is already $\Acal$-quasiconvex,
namely
\[
  F_\eps^\rel[u]=\begin{cases}
    \displaystyle\int_{\Omega_1} \Qcal_{\Acal}f(x,u(x))\dd{x} & \text{if } u\in \Ucal_\eps,\\
    +\infty & \text{otherwise.}
  \end{cases}
\]
Then we construct for every $u\in \Ucal_0$ a (strongly convergent) thin-film recovery sequence $(u_j)_j\in \Ucal_{\eps_j}$ 
for given $\eps_j \todown 0$. 
This is accomplished by passing to a Fourier point of view and using a convergence of symbols in order to infer strong convergence of 
the corresponding projections. 

\begin{proposition} \label{prop:rec_sec}
Suppose that Assumptions~\assref{ass:crp}, \assref{ass:Ad_triangle}, \assref{ass:approx_extension} and \assref{ass:antisymmetry} are satisfied. 
For $u \in \Ucal_0$ and $\eps_j\todown 0$ as $j\to \infty$, there exists a sequence $u_j \in \Ucal_{\eps_j}$ ($j \in \N$) 
with $u_j \to u$ strongly in $L^p(\Omega_1;\R^m)$ for $j\to \infty$.
\end{proposition}

For the proof of the proposition we will need the auxiliary symbol
\[
  \tilde{\Abb}_0(\xi) := \left( \left\{ \begin{aligned}
      &[A^{(d)}]^i \xi_d              && \text{if $[A^{(d)}]^i\xi_d \neq 0$,} \\
      &\sum_{k=1}^{d-1} [A^{(k)}]^i \xi_k && \text{if $[A^{(d)}]^i\xi_d = 0$,}
    \end{aligned} \right\} \right)^{i = 1,\ldots,l},\qquad \xi\in \R^d,
\]
which differs from $\Abb_0(\xi)$ only on the hyperplane where $\xi_d = 0$. 
For $\xi_d=0$ one has $\tilde{\Abb}_0(\xi)=\Abb'(\xi')$.
We also denote by $\tilde{\Pbb}_0(\xi)$ the corresponding projection onto $\ker \tilde{\Abb}_0(\xi)$.

\begin{remark}
In contrast to $\Abb_0(\xi)$, the ``symbol'' $\tilde{\Abb}_0(\xi)$ is not a polynomial and so it does not 
correspond to a constant-coefficient differential operator.
\end{remark}

We start with the following lemma.

\begin{lemma} \label{lem:P_eps_conv}
Under the Assumptions~\assref{ass:crp} and \assref{ass:Ad_triangle} and with $\eps_j\todown 0$ as $j\to \infty$, 
the symbols $\Abb_{\eps_j}$ ``converge'' for $j\to \infty$ to the symbol $\tilde{\Abb}_0$ in the sense that 
$\Pbb_{\eps_j}(\xi) \to \tilde{\Pbb}_0(\xi)$ for all $\xi \in \R^d$.
\end{lemma}
\begin{proof}
The case $\xi_d = 0$ is clear, and by the positive $0$-homogeneity of the projections we may assume 
$\xi \in \Sbb^{d-1}$ with $\xi_d \neq 0$. In view of \eqref{lem:step1} we know that 
$\ker \Abb_{\eps_j}(\xi)=\ker \Abb_{\eps_j}(\xi)_+$, if $\eps_j$ is small enough. The basic idea is that we can write the kernel of $\Abb_{\eps_j}(\xi)_+$ as
\begin{align}\label{kernelAeps}
  \ker \Abb_{\eps_j}(\xi)_+ = \set{ v \in \R^m }{ v\perp\frac{\Abb_{\eps_j}(\xi)^i}
    { \abs{\Abb_{\eps_j}(\xi)^i}}  \text{ for all $i = 1,\ldots,r$ }},
\end{align}
where $\Abb_{\eps_j}(\xi)^i$ denotes the $i$th row of $\Abb_{\eps_j}(\xi)$ (or $\Abb_{\eps_j}(\xi)_+$). 
Notice that for sufficiently small $\eps_j$ it holds that
$\abs{\Abb_{\eps_j}(\xi)^i} \neq 0$.

We then show that for $j\to \infty$, the kernel in \eqref{kernelAeps} converges (in the sense that the projections converge) 
to the kernel of $\tilde{\Abb}_0(\xi)_+$, which coincides with $\ker \tilde{\Abb}_0(\xi)$ by Lemma~\ref{lem:lindep}.
Since $\eps \Abb_\eps(\xi)_+ \to A_+^{(d)}\xi_d$ for $\eps\to 0$,
\begin{equation} \label{eq:ker_normals_conv}
  \lim_{j\to \infty} \, \frac{\Abb_{\eps_j}(\xi)^i}{\abs{\Abb_{\eps_j}(\xi)^i}}
  = \frac{[A^{(d)}]^i \xi_d}{\absn{[A^{(d)}]^i \xi_d}}  = \frac{\Abb_0(\xi)^i}{\abs{\Abb_0(\xi)^i}}
  = \frac{\tilde{\Abb}_0(\xi)^i}{\abs{\tilde{\Abb}_0(\xi)^i}}
\end{equation}
for all $i = 1,\ldots, r$.

In the spirit of the Gram--Schmidt orthogonalization procedure, we inductively define for $i = 1,\ldots,r$ the set of vectors 
\begin{align*}
  v_j^i &:= \frac{\Abb_{\eps_j}(\xi)^i}{\abs{\Abb_{\eps_j}(\xi)^i}}
    - \sum_{k=1}^{i-1} \biggl( v_j^k \cdot
    \frac{\Abb_{\eps_j}(\xi)^k}{\abs{\Abb_{\eps_j}(\xi)^k}} \biggr) v_j^k, \\
  w_j^i &:= \frac{v_j^i}{\absn{v_j^i}}.
\end{align*}
Notice that the $w_j^i$ are well-defined since the rows of $\Abb_{\eps_j}(\xi)_+$ are linearly independent, and therefore
the $v_j^i$ are never zero. Analogously, define the collections $v_0^i$, $w_0^i$ ($i = 1,\ldots,r$), but with $\Abb_{\eps_j}(\xi)^i$ 
replaced by $\tilde{\Abb}_0(\xi)^i = [A^{(d)}]^i \xi_d$.

By construction, the $w_j^i$ ($i = 1,\ldots,r$) for fixed $j\in \N$ form an orthonormal system (the same holds true for the $w_0^i$), 
and~\eqref{eq:ker_normals_conv} implies that as $j\to \infty$,
\begin{align}\label{convergence_wji}
  w_j^i \to w_0^i  \qquad\text{for all $i = 1,\ldots,r$.}
\end{align}
We may write the projections $\Pbb_{\eps_j}(\xi)$ and $\tilde{\Pbb}_0(\xi)$ for any $v \in \R^m$ as
\[
  \Pbb_{\eps_j}(\xi)v = v - \sum_{k=1}^r (v \cdot w_j^i) w_j^i,
  \qquad
  \tilde{\Pbb}_0(\xi)v = v - \sum_{k=1}^r (v \cdot w_0^i) w_0^i,
\]
and use \eqref{convergence_wji} to conclude that $\Pbb_{\eps_j}(\xi)v \to \tilde{\Pbb}_0(\xi)v$ for $j\to \infty$.
\end{proof}

We can now turn to the proof of Proposition~\ref{prop:rec_sec}.

\begin{proof}[Proof of Proposition~\ref{prop:rec_sec}]
The proof is divided into several steps.

\proofstep{Step 1: Approximative extension to $\T^d$.}
Take $u \in \Ucal_0$. By Assumption~\assref{ass:approx_extension}, we find $v \in L^p(\T^d;\R^m)$ with $\Acal_0 v = 0$ in $\T^d$ 
and $v$ close to $u$ in the $L^p(\Omega_1;\R^m)$-norm (recall that without loss of generality we assume $\omega \subset\subset Q^{d-1}$). 
If we can show the assertion for $v$, then we may conclude the statement of the theorem by a diagonal argument.

\proofstep{Step 2: Splitting of $v$.}
First additionally assume that $v$ is smooth. The $\Acal_0$-freeness of $v$ implies
\[
  \Abb_0(\xi) \hat{v}(\xi) = 0  \qquad \text{for all $\xi \in \Z^d$.}
\]
Now split $\hat{v}$ into
\[
  \hat{v}^{(1)}(\xi) := \begin{cases}
                      \hat{v}(\xi)  & \text{if $\xi_d \neq 0$,} \\
                      0             & \text{if $\xi_d = 0$,}
                    \end{cases}
  \qquad\text{and}\qquad
  \hat{v}^{(2)}(\xi) := \begin{cases}
                      0             & \text{if $\xi_d \neq 0$,} \\
                      \hat{v}(\xi)  & \text{if $\xi_d = 0$,}
                    \end{cases}
\]
so that
\[
  v = v^{(1)} + v^{(2)}.
\]
The functions $v^{(1)}, v^{(2)}$ are still smooth and satisfy $\Acal_0 v^{(1)} = \Acal_0 v^{(2)} = 0$ in $\T^d$. 
The smoothness follows for example by observing that
\begin{equation*} 
  v^{(2)}(x) = v^{(2)}(x') = \int_0^1 v(x',s) \dd s,
  \qquad x  \in Q^d,
\end{equation*}
which can be proved by comparing Fourier coefficients.

We will now show the existence of sequences $(v_j^{(1)})_{j}$ and $(v_j^{(2)})_{j}$ with $v_j^{(1)} \to v^{(1)}$, $v_j^{(2)} \to v^{(2)}$ 
strongly in $L^p(Q^d;\R^m)$ as $j\to \infty$ and $\Acal_{\eps_j} v_j^{(1)} = 0$ in $\T^d$, $\Acal_{\eps_j} v_j^{(2)} = 0$ in $Q^d$.

\proofstep{Step 3: The part $v^{(1)}$.}
For $v^{(1)}$ we have
\[
  \tilde{\Abb}_0(\xi) \hat{v}^{(1)}(\xi) = 0  \qquad \text{for all $\xi \in \Z^d$.}
\]
Set
\[
  v_j^{(1)} := \Pcal_{\eps_j} v^{(1)},  \qquad\text{i.e.}\qquad
  \hat{v}_j^{(1)}(\xi) = \Pbb_{\eps_j}(\xi) \hat{v}^{(1)}(\xi),
\]
where $\Pcal_{\eps_j}$ is the projection onto the kernel of $\Acal_{\eps_j}$ as in Theorem~\ref{theorem_projection}; 
we let $\Pbb_{\eps_j}(0) = I_m$ with $I_m$ the identity map on $\R^m$.

Recall that for smooth functions the Fourier inversion formula holds, so that
\[
  v_j^{(1)}(x) = \sum_{\xi \in \Z^d} \hat{v}_j^{(1)}(\xi) \ee^{2\pi\ii x\cdot\xi}
  = \sum_{\xi \in \Z^d} \Pbb_{\eps_j}(\xi) \hat{v}^{(1)}(\xi) \ee^{2\pi\ii x\cdot\xi},
  \qquad x \in Q^d.
\]
Thus, 
\begin{align}\label{estimate1}
  \normb{v_j^{(1)} - v^{(1)}}_{L^p(Q^d;\R^m)}^p 
&= \int_{Q^d} \absB{ \sum_{\xi \in \Z^d \setminus \{0\}} \bigl(\Pbb_{\eps_j}(\xi) \hat{v}^{(1)}(\xi) 
- \hat{v}^{(1)}(\xi)\bigr) \ee^{2\pi\ii x\cdot\xi}\, }^p \dd x\nonumber\\
  &\leq \Big[ \sum_{\xi \in \Z^d \setminus \{0\}} \absb{\Pbb_{\eps_j}(\xi)\hat{v}^{(1)}(\xi) - \hat{v}^{(1)}(\xi)}  \Big]^p.
\end{align}
The individual terms of the last sum converge to zero, because
\[
  \Pbb_{\eps_j}(\xi)\hat{v}^{(1)}(\xi) \to \tilde{\Pbb}_0(\xi)\hat{v}^{(1)}(\xi) = \hat{v}^{(1)}(\xi)  \qquad\text{as $j\to \infty$,}
\]
by Lemma~\ref{lem:P_eps_conv}. Then, the rapid decay of Fourier coefficients of smooth functions together with the uniform boundedness
 of the projection matrices $\Pbb_{\eps_j}(\xi)$ allows us to invoke a series version of Lebesgue's dominated convergence theorem to 
see that the right-hand side in $\eqref{estimate1}$ converges to zero as $j \to \infty$. Hence we may conclude that 
$v_j^{(1)} \to v^{(1)}$ strongly in $L^p(Q^d;\R^m)$.

\proofstep{Step 4: The part $v^{(2)}$.}
Turning to $v^{(2)}$, we recall that by construction $v^{(2)}(x) = v^{(2)}(x')$.  

We now construct a matrix $(A^{(d)})^\dagger \in \R^{m \times l}$, which in fact is the Moore--Penrose pseudoinverse of $A^{(d)}$, such that
\begin{align}\label{Moore_Penrose}
  A^{(d)} (A^{(d)})^\dagger = \left[ \begin{array}{c|c}
    I_r & \; 0 \\ \hline
    0   & \; 0
  \end{array}\right]
  \qquad \text{$\in \R^{l \times l}$,}
\end{align}
where $I_r$ is the identity matrix in $\R^{r \times r}$. 
Indeed, we observe that $A^{(d)}_+\in \R^{r\times m}$ has full row-rank $r$ and is hence surjective as a mapping into $\R^r$. Therefore, one can find 
$(A^{(d)}_+)^\dagger\in \R^{m\times r}$ such that
$A_+^{(d)}(A_+^{(d)})^\dagger=I_r$.
If we set
\[
  (A^{(d)})^\dagger = \left[ \begin{array}{c|c}
   (A^{(d)}_+)^\dagger & \: 0
  \end{array}\right]
  \qquad \text{$\in \R^{m \times l}$,}
\]
we get $\eqref{Moore_Penrose}$.

With $(A^{(d)})^\dagger$ at hand, set
\[
  v_j^{(2)}(x) := v^{(2)}(x') - {\eps_j} x_d (A^{(d)})^\dagger \Acal' v^{(2)}(x'), \qquad x\in Q^d.
\]
Notice that $v^{(2)}_j\notin C^\infty(\T^d;\R^m)$; by definition $v^{(2)}_j$ has a jump in $x_d$-direction over the gluing boundary. 
However, $v^{(2)}_j\in C^\infty(Q^d;\R^m)$. In view of $A^{(d)}_- = 0$ and the fact that $v^{(2)}$ is $\Acal_0$-free, we find for $x'\in Q^{d-1}$ that
\[
  \Acal'_- v^{(2)}(x')  = (\Acal_0)_- v^{(2)}(x') = 0.
\]
Thus,
\[
  A^{(d)} (A^{(d)})^\dagger \Acal' v^{(2)}(x') = \left[ \begin{array}{c|c}
    I_r & \; 0 \\ \hline
    0   & \; 0
  \end{array}\right]
  \cdot
  \left[ \begin{array}{c}
    \Acal'_+ v^{(2)}(x') \\ \hline
    0
  \end{array}\right]
  =
  \Acal' v^{(2)}(x').
\]
If we apply $\Acal_{\eps_j}$ to $v_j^{(2)}$, we get for all $x\in Q^d$ that
\begin{align*}
  \Acal_{\eps_j} v^{(2)}_j(x) &= \Acal' v^{(2)}_j(x)
    + \frac{1}{\eps_j} A^{(d)} \partial_d v^{(2)}_j(x) \\
  &= \Acal' v^{(2)}(x')
    - {\eps_j} x_d \Acal'\bigl[ (A^{(d)})^\dagger
    \Acal' v^{(2)}\bigr](x')
    - A^{(d)} (A^{(d)})^\dagger \Acal' v^{(2)}(x') \\
  &=\Acal' v^{(2)}(x') - {\eps_j} x_d \Acal'
    \bigl[ (A^{(d)})^\dagger \Acal' v^{(2)}\bigr](x')  - \Acal' v^{(2)}(x') \\
  &= - \eps_j x_d \Acal'\bigl[ (A^{(d)})^\dagger
    \Acal' v^{(2)} \bigr](x').
\end{align*}
Owing to Assumption~\assref{ass:antisymmetry} and the symmetry of second derivatives,
\begin{align*}
 \Acal'\bigl[ (A^{(d)})^\dagger \Acal' v^{(2)}\bigr]= \sum_{k,l=1}^{d-1} A^{(k)}(A^{(d)})^\dagger A^{(l)}\partial_k\partial_jv^{(2)}=0.
\end{align*}
So, 
\begin{align*}
 \Acal_{\eps_j}v^{(2)}_j = 0\qquad \text{in $Q^d$ (classically and distributionally).}
\end{align*}
Combining this with the observation that $v_j^{(2)} \to v^{(2)}$ in $L^p(Q^d;\R^m)$ gives the sought recovery sequence for $v^{(2)}$.

\proofstep{Step 5: $v$ non-smooth.}
If $v$ of Step~1 is not smooth, we take a sequence of smooth mollifying kernels $(\eta_n)_n \subset C_c^\infty(\R^d)$ 
with $\eta_n \to \delta_0$ as $n \to \infty$ and consider the functions
\[
  v_n(x) := (\eta_n \conv v)(x) = \int_{\R^d} v(x-y) \eta_n(y) \dd y,
  \qquad x \in \R^d,
\]
where $v$ here is considered as a $Q^d$-periodic function on all of $\R^d$. All the $v_n$ are smooth, $Q^d$-periodic and 
still $\Acal_0$-free on $\T^d$, because the operations of convolution and taking derivatives commute. Hence, the above 
reasoning applies and we find a recovery sequence for each $v_n$. Since also $v_n \to v$ strongly in $L^p(Q^d;\R^m)$, the 
existence of a recovery sequence for $v$ follows by a diagonal argument. By the reasoning in Step~1 this concludes the proof 
of the proposition.
\end{proof}

With Proposition \ref{prop:rec_sec} we finally obtain the upper bound.


\begin{proof}[Proof of Theorem~\ref{mainresult_Gamma}~(ii)]
The proof follows by joining Proposition \ref{prop:rec_sec} and the relaxation result of Theorem 1.1 in \cite{BFL00} 
through a diagonalization procedure. 

More precisely, by Proposition \ref{prop:rec_sec} we find for every $u\in \Ucal_0$ a sequence 
$u_j \in \Ucal_{\eps_j}$ ($j \in \N$) such that $u_j \to u$ in $L^p(\Omega_1;\R^m)$. 
Since $-\Qcal_\Acal f$ is a normal integrand satisfying the growth condition
\begin{align*}
 -\Qcal_\Acal f(x, v)\geq -C(1+\abs{v}^p)\qquad\text{ for almost every $x\in \Omega_1$ and all $v\in \R^m$},
\end{align*}
which is an immediate consequence of $\eqref{f_growth/coercivity}$,
and owing to the lower semicontinuity of $-\Qcal_\Acal f(x,\frarg)$ 
for almost all $x \in \Omega_1$,
we may infer that the functional $v\mapsto \int_{\Omega_1} \Qcal_\Acal f(x, v(x))\dd{x}$ is upper semicontinuous with respect 
to strong convergence in $L^p(\Omega_1;\R^m)$, see for instance Theorem~6.49 in~\cite{FL07}.
Hence,
\[
  \limsup_{j\to\infty} \int_{\Omega_1} \Qcal_{\Acal} f\bigl(x,u_j(x)\bigr) \dd{x} \leq \int_{\Omega_1} \Qcal_{\Acal} f\bigl(x,u(x)\bigr) \dd{x}.
\]
In view of the relaxation result in~\cite{BFL00}, for each $j\in \N$ there exists a sequence 
$(u_{j}^k)_{k}\subset \Ucal_{\eps_j}$ with $u_j^k \weakly u_j$ as $k\to \infty$ such that
\begin{align*}
\lim_{k\to \infty} F_{\eps_j}[u_{j}^k] & =\lim_{k\to \infty}\int_{\Omega_1}  f\bigl(x,u_{j}^{k}(x)\bigr)\dd{x} 
= \int_{\Omega_1} \Qcal_{\Acal_{\eps_j}} f\bigl(x,u_j(x)\bigr)\dd{x}\\ 
& =\int_{\Omega_1} \Qcal_\Acal f\bigl(x,u_j(x)\bigr)\dd{x}.
\end{align*}
Here we used Lemma~\ref{lem:Aeps-Adelta_quasiconvexity}, which implies that $\Qcal_{\Acal_{\eps_j}}f = \Qcal_{\Acal}f$ for all $j\in \N$.
Finally, we pick appropriate $k(j)$ for $j \in \N$ to conclude with the sequence $(u_j^{k(j)})_j$.
\end{proof}

\section{A-priori locality of the lower bound implies $\Gamma$-convergence}
\label{sc:locality}

The question of locality of the $\Gamma$-limit $F_0$ of $F_\eps$ is an interesting open problem that was already pointed out~in \cite{BFM03} for 
the gradient case.
If one however a-priori assumes that the $\Gamma$-limit is local (along with some technical requirements), we can prove that the integrand in $F_0$ 
can be fully identified and turns out to be equal to the upper bound $\Qcal_\Acal f$.

In all of the following, we assume that $f$ is independent of $x$ for simplicity. Consider $f:\R^m\to \R$ and assume that $F_0=\Gamma$-$\lim_{\eps\to 0} F_\eps$ exists and is local (see below). Then as a consequence of Proposition~\ref{prop:locality} below we find that
\begin{align*}
F_0[u]=\begin{cases} \displaystyle\int_{\Omega_1} \Qcal_\Acal f(u(x))\dd{x}, & u\in \Ucal_0, \\ \infty, & \text{otherwise.} \end{cases}
\end{align*}


To show this, let us define for any open subset $D$ of $\Omega_1$ and $u\in \Ucal_0$,
\begin{align*}
F_0^-[u;D]:=\inf\, \biggl\{\liminf_{j\to 0} \int_{D} f(u_j(x))\dd{x} : u_j \in \Ucal_{\eps_j}^D, u_j\weakly u \text{ in $L^p(D;\R^m)$}\biggr\}
\end{align*}
with $\Ucal_\eps^D:=\{u\in L^p(D;\R^m) : \Acal_{\eps} u=0 \text{ in $D$}\}$ for $\eps>0$.
\begin{proposition}\label{prop:locality}
Suppose $\Acal$ and $f:\R^m\to \R$ satisfy the conditions of Theorem~\ref{mainresult_Gamma} and $\Qcal_\Acal f$ is continuous. Let $u\in \Ucal_0$. 
If $F_0^-$ is a local integral functional, that is, $F_0^-[u;\frarg]$ is equal to a Radon measure absolutely continuous with 
respect to the Lebesgue measure restricted 
to $\Omega_1$, or, more precisely, there exists a function $g\in L^1(\Omega_1)$
such that $F_0^-[u;\frarg] = g\,\Lcal^d\lfloor \Omega_1$, then
\begin{align*}
 g(x)= \Qcal_\Acal f(u(x))
\end{align*}
for almost every~$x\in \Omega_1$.
\end{proposition}

\begin{proof}
It is enough to prove $g\geq \Qcal_\Acal f(u)$. The other inequality follows directly from Theorem~\ref{mainresult_Gamma} (ii).
Let $x_0\in \Omega_1$ be a $p$-Lebesgue point of $u$, i.e.\
\begin{align}\label{pLebesguepoint}
 \lim_{r\to 0^+} \frac{1}{r^n} \int_{Q(x_0, r)} \abs{u(x)-u(x_0)}^p \dd{x}=0,
\end{align}
where $Q(x_0,r) := x_0 + (-r,r)^d$. Further, in view of the Besicovitch derivation theorem, assume
\begin{align*}
 g(x_0)=\lim_{r\to 0^+} \frac{F_0^-(u;Q(x_0,r))}{r^n}<\infty
\end{align*} 
and suppose that $F_0^-(u;\partial Q(x_0,r))=0$ for the chosen radii $r>0$.
For fixed $r$, we consider $u_j^r \in U_{\eps_j}^{Q(x_0,r)}$ $(j\in \N)$ such that $u_j^r\weakly u$ in $L^p(Q(x_0,r);\R^m)$ as $j\to\infty$ and
\begin{align*}
\lim_{j\to \infty} \int_{Q(x_0,r)} f(u_j^r(x)) \dd{x} \leq F_0^-[u;Q(x_0,r)] + r^{n+1}.
\end{align*}
Then,
\begin{align*}
g(x_0) \geq \liminf_{r\to 0^+} \lim_{j\to \infty} \frac{1}{r^n} \int_{Q(x_0,r)} f(u_j^r(x))\dd{x}
=\liminf_{r\to 0^+} \lim_{j\to \infty}\int_{Q(0,1)} f(u(x_0) + v_j^r(y))\dd{y}.
\end{align*}
In the last equality we performed the change of variables $y:=(x-x_0)/r$ and set $v_j^r(y)=u_j^r(x_0+ry)-u(x_0)$ for $y\in Q(0,1)$. 
Then $v_j^r\weakly 0$ in $L^p(Q(0,1);\R^m)$ as $j\to \infty$, since for every $w\in L^{p'}(Q(0,1);\R^m)$ one obtains that
\begin{align*}
&\absBB{\int_{Q(0,1)} v_j^r(y) w(y)\dd{y}}\\ &\quad\leq \absBB{\int_{Q(0,1)} (u_j^r(x_0+ry)-u(x_0+ry))w(y)\dd{y}} + \absBB{\int_{Q(0,1)} (u(x_0+ry)-u(x_0))w(y)\dd{y}} \\
&\quad\leq \absBB{\frac{1}{r^n} \int_{Q(x_0,r)} (u_j^r(x)-u(x))w((x-x_0)/r)\dd{x}} \\ 
& \qquad\qquad\qquad\qquad\qquad+ \norm{w}_{L^{p'}(Q(0,1);\R^m)} \biggl(\frac{1}{r^n} \int_{Q(x_0,r)}\abs{u(x)-u(x_0)}^p\dd{x}\biggr)^{1/p}.
\end{align*}
Due to the weak convergence of $(u_j^r)_j$ to $u$ in $L^p(Q(x_0,r);\R^m)$ and \eqref{pLebesguepoint}, the right-hand side in the above expression tends 
to zero for $j\to \infty$.

By a diagonalization argument we may now extract a sequence $\hat{v}_k \in U_{\Acal_{\eps_{j(k)}}}$ ($k \in \N$) such that $\hat{v}_k\weakly 0$ in $L^p(Q(0,1);\R^m)$ and 
\begin{align*}
g(x_0)\geq \lim_{k\to \infty} \int_{Q(0,1)} f(u(x_0) + \hat{v}_k(y))\dd{y}.
\end{align*}
Hence, by the definition of $\Acal_{\eps_{j(k)}}$-quasiconvexity, and with $\bar{v}_k := \int_{Q(0,1)}\hat{v}_k\dd{y}$,
\begin{align*}
g(x_0) &\geq \liminf_{k\to \infty} \int_{Q(0,1)} \Qcal_{\Acal_{\eps_{j(k)}}} f\bigl(u(x_0) + \bar{v}_k +\hat{v}_k(y) - \bar{v}_k \bigr)\dd{y}\\ 
&\geq \liminf_{k\to \infty} \Qcal_{\Acal_{\eps_{j(k)}}} f\bigl(u(x_0) + \bar{v}_k \bigr) \\ &
=\liminf_{k\to \infty} \Qcal_{\Acal} f\bigl(u(x_0) + \bar{v}_k \bigr)  =\Qcal_\Acal f(u(x_0)).
\end{align*}
Here we used that $\Qcal_{\Acal_\eps}f=\Qcal_\Acal f$ for all $\eps>0$, and $\bar{v}_k\to 0$ as $k\to \infty$ together with the continuity of $\Qcal_\Acal f$. This finishes the proof.
\end{proof}

\begin{remark}
Notice that requiring continuity of $\Qcal_\Acal f$ in Proposition~\ref{prop:locality} is not restrictive for $\Acal=\diverg$ and $\Acal=\curl$ or if $f$ is $\Acal$-quasiconvex.
\end{remark}

\section{Application: Thin films in nonlinear elasticity}\label{sec:application}

The energies governing (hyper-)elastic bulk bodies take the form of integral functionals depending on deformation gradients.
Since gradients are essentially the $\curl$-free vector fields, there is an alternative way of modeling an elastic energy by imposing a suitable PDE constraint. 
In what follows we compare these two modeling strategies when passing to the thin-film limit.
It turns out that the $\curl$-free formulation has the advantage of supplying strictly more information 
and is in fact equivalent to models accounting for bending through a so-called Cosserat vector~\cite{BFM03, BFM09}.

In this section we work within a three-dimensional setting and assume that $\Omega_\eps=\omega\times (0,\eps)\subset \R^3$ models 
the reference configuration of a film of thickness $\eps>0$ where the cross section $\omega\subset\R^2$ is a bounded and simply 
connected Lipschitz domain.
The elastic energy density $f:\R^{3\times 3}\to \R$ is supposed to be continuous and to satisfy $p$-growth and $p$-coercivity in the 
sense of~\eqref{f_growth/coercivity}.
For simplicity we dispense with the explicit dependence of $f$ on the space variable, which corresponds to assuming a
homogeneous material response. A practical example of such an energy density is $f(M)=\dist^p(M, \text{\rm SO}(3))$ with 
$M\in \R^{3\times 3}$ and $\text{\rm SO}(3)$ the rotation group. 
For $p=2$ this function $f$ has quadratic growth and meets the usual assumptions in 
geometrically nonlinear elasticity. For example, with the above choice $f$ is frame indifferent.

A classical starting point for deriving lower dimensional membrane models from $3$d elasticity~\cite{LR95, LR00, LR96, Bra02} is the energy functional given by
\begin{align*}
G_\eps^{\rm cl}[v]
= \frac{1}{\eps}\int_{\Omega_\eps} f\bigl(\nabla v(y)\bigr)\dd{y}, \qquad v\in W^{1,p}(\Omega_\eps;\R^3).
\end{align*}
On the other hand, in the approach based on $\curl$-free vector fields we seek to investigate
\begin{align}\label{tilde_Feps}
G_\eps[K]=\begin{cases}
\displaystyle\frac{1}{\eps}\int_{\Omega_\eps} f(K(y))\dd{y} &\text{if $\curl K=0$  
in $\Omega_\eps$,}\\ +\infty& \text{otherwise,}\end{cases}
\end{align} 
with $K\in L^p(\Omega_\eps;\R^{3\times 3})$.
After the thin-film rescaling $\eqref{rescaling_variables}$, which allows us to work on the fixed domain $\Omega_1$, 
the energy $G_\eps^{\rm cl}$ turns into
\begin{align}\label{overline_Geps}
F_\eps^{\rm cl}[u]=\int_{\Omega_1} f\Bigl(\nabla' u(x)\,\Big|\, 
\frac{1}{\eps} \,\partial_3 u(x) \Bigr)\dd{x}, \qquad u\in W^{1,p}(\Omega_1;\R^3),
\end{align} 
where we use the notation $\nabla'u=(\partial_1 u\,|\,\partial_2 u)$,
%
while in $\eqref{tilde_Feps}$ the change of variables $\eqref{rescaling_variables}$ together with $H(x)=K(y)=K(x', \eps x_d)$ for $x\in \Omega_1$ implies
\begin{align*}
F_\eps[H]=\begin{cases}
\displaystyle\int_{\Omega_1} f\bigl(H(x)\bigr)\dd{x} &\text{if $\curl_\eps H=0$  in $\Omega_1$,} \\
+\infty &\text{otherwise},\end{cases}\qquad H\in L^p(\Omega_1;\R^{3\times 3}).
\end{align*} 

%

The reduced limit functional $F_0^{\rm cl}$, 
obtained by $\Gamma$-convergence of the family $F_\eps^{\rm cl}$ regarding the weak $W^{1,p}$-topology, is given by $F_0^{\rm cl}:W^{1,p}(\omega;\R^{3})\to \R$ with
\begin{align*}
F_0^{\rm cl}[u]=\displaystyle\int_{\omega} \Qcal_2\Bigl(\min_{b\in \R^3} f\bigl(\partial_1 u(x')\big|\partial_2 u(x')\big|b\bigr)\Bigr)\dd{x'}, \qquad u\in W^{1,p}(\omega;\R^3).
\end{align*}
Here $\Qcal_2h$ stands for the operation of quasiconvexification in two dimensions of a function $h:\R^2\to \R$. For the proofs and 
further details we refer to \cite{LR95,Bra02}. 

Accounting for Example~\ref{ex:curl} (with $d=3$, $n=1$ and $\curl$ interpreted as $\nabla \times$) and under
the additional assumption that $f$ is asymptotically $\curl_0$-quasiconvex, i.e.\ $f=\Qcal^\infty_{\curl_0} f$,
the $\Gamma$-limit $F_0$ of $F_\eps$ regarding weak convergence in $L^p$ can be computed by Theorem~\ref{mainresult_Gamma} and Remark~\ref{rem:envelopes},
whereby we find
\begin{align}\label{F0}
F_0[H]=\begin{cases}\displaystyle\int_{\Omega_1} f\bigl(H(x)\bigr)\dd{x} &\text{if } H\in \Ucal_{F_0},\\
+\infty &\text{otherwise,}\end{cases}
\end{align} 
with 
\begin{align*}
\Ucal_{F_0} &:= \set{H\in L^p(\Omega_1;\R^{3\times 3})}{\curl_0 H=0 \text{ in $\Omega_1$}}\\
&=\set{H\in L^p(\Omega_1;\R^{3\times 3})}{H=(\nabla'h|H_3),\, h\in W^{1,p}(\omega;\R^3),\, H_3\in L^p(\Omega_1;\R^3)}.
\end{align*}
Notice that in contrast to $F_0^{\rm cl}$, the $\Gamma$-limit $F_0$ is substantially three-dimensional, and therefore contains 
strictly more information than the purely two-dimensional $F_0^{\rm cl}$.

Summarizing, the essential drawback of the classical approach is that the 
information about the weak limit of $\frac{1}{\eps}\partial_3 u_\eps$ is lost as $\eps\to 0$.
For this reason the authors of~\cite{BFM09, BFM03} provide an extended analysis of the gradient formulation that keeps 
track of exactly this quantity, which is called the Cosserat vector and reveals additional insight into bending behavior
beyond mere dimension reduction. In fact, Bouchitt\'e, Fonseca and Mascarenhas~\cite{BFM09} study the (already rescaled) functional 
$I_\eps[u,b]\colon W^{1,p}(\Omega_1;\R^3)\times L^p(\Omega_1;\R^3)\to \overline{\R}$ defined as
\begin{align*}
I_\eps[u,b]=\begin{cases}
\displaystyle\int_{\Omega_1} f\Bigl(\nabla' u(x)\,\Big|\,\frac{1}{\eps}\partial_3 u(x) \Bigr)\dd{x} 
&\text{if } b= \frac{1}{\eps}\partial_3 u, \\
+\infty &\text{otherwise,}
\end{cases}
\end{align*}
and give a representation of the $\Gamma$-limit $I_0$ as $\eps$ tends to $0$. They conjectured that $I_0$ is nonlocal in general, 
but were able to find an integral representation for the special class of cross-quasiconvex integrands.
We point out that cross-quasiconvexity first appeared in the literature under the 
name ``joint quasiconvexity/convexity'',~see~\cite{FKP94}. Here we give only the definition and refer to~\cite{LR00, BFM09, BFM03} for further details.

\begin{definition}\label{def:cross-qc}
A function $\bar{f}:\R^{m\times{(d-1)}}\times\R^m\to \R$ is called cross-quasiconvex, if 
\begin{align}\label{estimate_cross-qc}
\dashint_{Q^{d-1}} \bar{f}(V+\nabla w(x'), v + z(x'))\dd{x'}\geq \bar{f}(V,v)
\end{align}
for all $(V, v)\in \R^{m\times{(d-1)}}\times \R^m$ and all $w\in C^\infty(\T^{d-1};\R^m)$ and $z\in C^\infty(\T^{d-1};\R^m)$ 
with $\int_{Q^{d-1}}z\dd{x'}=0$.
\end{definition}

\begin{remark}
If $\bar{f}$ is a continuous function satisfying the growth condition $\abs{\bar{f}(V,v)}\leq C (1 + \abs{(V,v)}^p)$ 
for all $(V, v)\in \R^{m\times{(d-1)}}\times \R^m$, then by a density 
argument~\eqref{estimate_cross-qc} holds for all $w\in W^{1,p}(\T^{d-1};\R^m)$ and all $z\in L^p(Q^{d-1};\R^m)$ with zero mean value.  
\end{remark}

The next lemma establishes the relation between cross- and $\curl_0$-quasiconvexity. 
 
\begin{lemma}\label{lem:crossqc}
A function $f:\R^{3\times 3}\to \R$ is $\curl_0$-quasiconvex if and only if $\bar{f}:\R^{3\times 2}\times \R^3\to \R$ defined 
through $\bar{f}(V,v)=f(V|v)$ with $V\in \R^{3\times 2}$ and $v\in \R^3$ is cross-quasiconvex.
\end{lemma}

\begin{proof}
First, assume that $\bar{f}$ is cross-quasiconvex. Let $(V|v)\in \R^{3\times 3}$ and 
$W\in C^\infty(\T^3;\R^{3\times 3})\cap \ker_{\T^3}\curl_0$ with $\int_{Q^3}W\dd{x}=0$. 
Then, $W$ can be respresented as $W=(\nabla w | z)$ with $w\in C^\infty(\T^2;\R^3)$ and $z\in C^\infty(\T^3;\R^3)$ such that $\int_{Q^3}z\dd{x}=0$, 
and
\begin{align*}
 \dashint_{Q^3} f\bigl((V|v)+ W(x)\bigr)\dd{x}= \dashint_{Q^3}\bar{f}\bigl(V+ \nabla w(x'), v+z(x)\bigr)\dd{x}.
\end{align*}
Since $\bar{f}$ is convex in its second argument by Proposition~4.4 of~\cite{LR00}, we may apply Jensen's inequality to infer
 \begin{align*}
 \dashint_{Q^3} f\bigl((V|v)+ W(x)\bigr)\dd{x} &\geq \dashint_{Q^2}\bar{f}\Bigl(V+ \nabla w(x'), v + \int_0^1 z(x', x_3)\dd{x_3}\Bigr)\dd{x'}\\
&\geq \bar{f}(V, v)= f(V|v).
\end{align*}
For the last inequality we used that the mapping $x'\to\int_0^1 z(x',x_3)\dd{x_3}$ lies in $C^{\infty}(\T^2;\R^3)$ with mean value zero 
and exploited the cross-quasiconvexity of $\bar{f}$.

If $f$ is $\curl_0$-quasiconvex we find for $w\in C^\infty(\T^{2};\R^m)$ and $z\in C^\infty(\T^{2};\R^m)$ 
with $\int_{Q^{2}}z\dd{x}=0$ (both identified with their constant extensions in the $x_3$-variable) that $\curl_0 (\nabla w|z) =0$ in $\T^3$ and 
$\int_{Q^3}(\nabla w| z)\dd{x}=0$. Hence, for any 
$(V,v)\in \R^{3\times 2}\times \R^3$,
 \begin{align*}
 \dashint_{Q^2} \bar{f}\bigl(V+ \nabla w(x'), v+ z(x') \bigr)\dd{x'}
= \dashint_{Q^3}f\bigl((V| v)  + (\nabla w(x)| z(x))\bigr)\dd{x} \geq f(V|v) = \bar{f}(V,v),
\end{align*}
which shows the cross-quasiconvexity of $\bar{f}$.
\end{proof}

Observing that
\begin{equation}\label{equivalence_Ieps-Feps}
I_\eps[u,b]=F_\eps[H] \qquad\text{for}\qquad H=[\nabla'u \,|\,b] \in L^p(\Omega_1;\R^{3\times 3}),
\end{equation}
we find that $F_\eps$ represents an equivalent formulation for the elastic energy in the model with bending moment. 
Then, with regard to $\eqref{equivalence_Ieps-Feps}$ and $\eqref{F0}$, Theorem~\ref{mainresult_Gamma} has the 
following implication (compare Theorem~2.2/Pro\-position~2.4 of \cite{BFM09}).

\begin{corollary}\label{cor:elastic_thinfilm_limit}
Let $f$ be asymptotically $\curl_0$-quasiconvex, i.e.\ $f=\Qcal^\infty_{\curl_0} f$, and satisfy \eqref{f_growth/coercivity}.
Then, the $\Gamma$-limit with respect to weak convergence in $W^{1,p}(\Omega_1;\R^3)\times L^p(\Omega_1;\R^3)$ of the 
functional $I_\eps$ as $\eps\to 0$ is given
by $I_0\colon W^{1,p}(\omega;\R^3)\times L^p(\Omega_1;\R^3)\to \R$ with
\begin{align*}
I_0[u,b]= \int_{\Omega_1}  f\bigl(\nabla' u(x')\,\big|\,b(x) \bigr)\dd{x}
=\int_0^1\int_\omega f\bigl(\nabla'u(x')\,\big|\,b(x',x_3) \bigr)\dd{x'}\dd{x_3}.
\end{align*}
\end{corollary}

Hence, Theorem~\ref{mainresult_Gamma} can be viewed as an extension of the results on thin films with Cosserat vector 
in the gradient setting to the context of problems on $\Acal$-free vector fields.

%
%
\section*{Acknowledgments}
We thank Irene Fonseca, who suggested this problem, for many useful discussions and remarks on the subject.
Moreover we appreciate her careful reading of a preliminary version of the manuscript. 
Thanks are also due to Helmut Abels for his insightful comments on Fourier multipliers and to Stefan Kr\"omer for helpful 
discussions regarding the lower bound and for pointing out a mistake in a previous version of the paper.

Part of this work was carried out during a visit of F.~R.\ to Carnegie Mellon University (CMU) and to 
Centro de Matem\'atica e Aplica\c{c}\~{o}es of Universidade Nova de Lisboa (CMA/UNL), and a visit of C.~K.\ 
to the Oxford Centre for Nonlinear PDE (OxPDE). C.~K.\ was supported by the 
Funda\c{c}\~{a}o para a Ci\^{e}ncia e a Tecnologia 
(FCT) through the ICTI CMU--Portugal Program in Applied Mathematics and UTA-CMU/MAT/0005/2009, and by OxPDE. F.~R.\ gratefully
acknowledges the support of OxPDE through the EPSRC Science and Innovation award to OxPDE (EP/E035027/1), and of the Center 
for Nonlinear Analysis at CMU.

%
%


\end{document}